\documentclass[12pt]{amsart}
\usepackage{amssymb}
\usepackage{amsfonts}

\newtheorem{cl}{Claim}
\newtheorem{theorem}{Theorem}[section]
\newtheorem{proposition}[theorem]{Proposition}
\newtheorem{coro}[theorem]{Corollary}
\newtheorem{lemma}[theorem]{Lemma}
\newtheorem{example}[theorem]{Example}

\newcommand\C{\mathbb{C}}
\newcommand\N{\mathbb{N}}
\newcommand\R{\mathbb{R}}

\newcommand{\supp}{\operatorname{supp}}

\begin{document}

\title[Gaussian heat kernel bounds via Phragm\'en-Lindel\"of theorem ]{%
Gaussian heat kernel upper bounds via Phragm\'en-Lindel\"of theorem }
\author{Thierry Coulhon}
\thanks{TC's research was partially supported by the European
Commission
(IHP Network ``Harmonic Analysis and Related Problems''
2002-2006, Contract HPRN-CT-2001-00273-HARP)}
\author{Adam Sikora}
\thanks{AS's research was partially supported by an 
Australian Research Council (ARC) Discovery Grant DP 0451016 and New 
Mexico State University Summer Research Award}
\address{Thierry Coulhon, D\'epartement de Math\'ematiques, Universit\'e de
Cergy-Pontoise, Site de Saint-Martin, 2, rue Adolphe Chauvin, F 95302
Cergy-Pontoise Cedex, FRANCE}
\email{Thierry.Coulhon@math.u-cergy.fr}
\address{Adam Sikora, Department of Mathematical Sciences, New Mexico State
University, Las Cruces, NM 88003-8001, USA}
\email{asikora@nmsu.edu}
\date{\today}

\begin{abstract}
We prove that in presence of $L^2$ Gaussian estimates, so-called
Davies-Gaffney estimates, on-diagonal upper bounds imply precise off-diagonal
Gaussian upper bounds for the kernels of analytic families of operators on metric
measure spaces.
\end{abstract}

\maketitle
\tableofcontents

\setcounter{section}{0} \renewcommand{\theenumi}{\alph{enumi}} %
\renewcommand{\labelenumi}{\textrm{(\theenumi)}} %
\numberwithin{equation}{section}

\noindent

\section{Introduction}

The study of properties of second-order self-adjoint differential operators
often depends on Gaussian upper bounds for the corresponding heat kernel, that is
the kernel of the semigroup generated by these operators. Gaussian heat
kernel upper bounds play a crucial role in the proofs of many results
concerning boundedness of the Riesz transform, convergence of Bochner-Riesz
means or boundedness of spectral multipliers, as well as problems related to
maximal regularity properties (see for example the articles  \cite{CD4, CD2, CD1, MR, ACDH, He, Ch, A1, A2,  CCO,  DOS, CD3} and the monograph \cite{O}). If $p_t(x,y)$ denotes the heat kernel
corresponding to a second-order differential elliptic or sub-elliptic
operator, then the typical Gaussian heat kernel upper bound is of the form 
\begin{equation}  \label{gauss}
0\le p_t(x,y)\le \frac{C}{V(x,\sqrt t)} \exp\left(-\frac{d^2(x,y)}{Ct}\right)
\end{equation}
for all $t>0$, $x,y$ ranging in the space where the operator acts. For instance,
if $p_t$ is the kernel corresponding to the Laplace-Beltrami operator on a
Riemannian manifold, then, at least in some favorable cases, one expects $V(x,r)$ to be equal to the volume of
the geodesic ball of radius $r$ and centered at $x$ and $d$ denotes the
Riemannian distance. In the standard approach, proofs of Gaussian estimates
are divided into two steps. First one obtains on-diagonal estimates 
\begin{equation}  \label{diag}
p_t(x,x) \le \frac{C}{V(x,\sqrt t)}
\end{equation}
for all $t,x$.
Then the theory says that one can automatically improve on-diagonal bounds
by adding the Gaussian factor $\exp\left(-\frac{d^2(x,y)}{Ct}\right) $ and
obtain this way Gaussian bounds (\ref{gauss}). There are basically three
known methods to derive Gaussian bounds from on-diagonal bounds (\ref{diag}%
): Davies's perturbation method (see \cite{Da}, \cite{DP}, \cite{Co1}), the
integrated maximum principle (see \cite{Gr0}, \cite{Gr1}, \cite{Gr}) and
finite propagation speed for the wave equation (see \cite{Si1, Si2}).

The main aim of the present paper is to introduce a new method for deducing Gaussian
bounds from uniform bounds (\ref{diag}), which relies mainly on the
Phragm\'en-Lindel\"of theorem. Our approach is closely related to the main
idea behind complex interpolation, and it shows that surprisingly the
Gaussian bounds and the complex interpolation results are of similar nature.
This allows us to look at the off-diagonal Gaussian bounds from a new
perspective. The use of Phragm\'en-Lindel\"of theorems for heat kernel
estimates was introduced to our knowledge in \cite{Da1}, see in particular
Lemma 9, see also \cite[Theorem~3.4.8, p.103]{Da}. In \cite{Da1}, Davies
uses Gaussian bounds for real time and the Phragm\'en-Lindel\"of technique
to obtain complex time estimates for the heat kernel, that is estimates for $%
p_z(x,y)$ for all $z\in \mathbb{C}_+$, where $\mathbb{C}_+$ is the complex
half-plane ${\mbox{\Small{\rm Re}}} z > 0$. Roughly speaking, in our
approach, we reverse the order of Davies's idea and we use the
Phragm\'en-Lindel\"of technique to obtain both real and complex time
Gaussian bounds. Our method yields more precise complex time Gaussian bounds
than in \cite{Da1} and \cite[Theorem~3.4.8, p.103]{Da}.

We express the assumptions of our results in terms of so-called
Davies-Gaffney estimates. To our knowledge, these estimates were formulated for the first time in \cite{Da2}, but according to  Davies himself the idea stems from \cite{Ga}. Davies-Gaffney
estimates hold for essentially all self-adjoint, elliptic or subelliptic
second-order differential operators including Laplace-Beltrami operators on
complete Riemannian manifolds, Schr\"odinger operators with real-valued
potentials and electromagnetic fields, and Hodge-Laplace operators acting on
differential forms (see Theorem \ref{negpot} and Section \ref{vb} below).
A discrete time version of the Davies-Gaffney estimate is
discussed in \cite{cgz}\footnote{%
However, it is not clear how to extend the methods of the present work to
the discrete time case, in order to replace the use of the rather technical
discrete integrated maximum principle as in \cite{cgz}.}. 
Davies-Gaffney estimates are also easy to obtain. For non-negative self-adjoint 
operators, they are equivalent with the finite speed
propagation property for the corresponding wave equation (see \cite{Si2} and
Section \ref{dg} below). We discuss this equivalence here as a simple but
illuminating application of the Phragm\'en-Lindel\"of technique.

Our approach allows us to obtain far reaching generalizations of the results
obtained in \cite{Da, DP, Co1, Gr0, Gr1, Gr, Si1, Si2}. In the present
paper, we do not have to assume anything about the nature of the
infinitesimal generator of the semigroup under consideration; in particular,
the generating operator does not have to be a second-order differential
operator, and the semigroup does not have to be Markov. Our method works
also for operators acting on differential forms and more generally on vector
bundles \footnote{%
It is worth noting at this point that there is a connection between
estimates of the heat kernel on 1-forms and the $L^p$ boundedness of the
Riesz transform for $p>2$ (see for example \cite{CD1, CD2, Si2}).}.
Actually, instead of considering the analytic semigroup $\{\exp(-zL)\colon \, z \in \mathbb{C}_+\}$ generated by some
non-negative self-adjoint operator $L$, we are
able to study any uniformly bounded analytic family of operators $%
\{\Psi(z)\colon \, z \in \mathbb{C}_+\} $. We do not have to assume that $%
\Psi$ has the semigroup property nor that $\Psi(z)$ is a linear operator, as
far as $\{\Psi(z)\colon \, z \in \mathbb{C}_+\} $ satisfies the
Davies-Gaffney estimates (see (\ref{DG2}) below). For example we can study
the  estimates for the gradient of the heat kernel in the same way as the
estimates for the heat kernel itself. We are also able to consider the
family given by the formula $\Psi(z)=\exp(-zL)-\exp(-z{L^0})$ where $L$, ${%
L^0}$ are different generators of analytic semigroups. For example, one can
consider the situation where $L$ is an operator with periodic coefficients
in divergence form and $L^0$ is its homogenization, to obtain Gaussian
estimates for the difference of the corresponding heat kernels $|p_t(x,y)-
p^0_t(x,y)|$.

Next, our methods have various applications in the theory of $L^p$ to $L^q$
Gaussian estimates, developed by Blunck and Kunstmann in \cite{B1, B2, BK1,
BK2, BK3}, see also \cite{LSV} and \cite{AM}. Blunck and Kunstmann call such estimates
generalized Gaussian estimates. They use generalized Gaussian estimates to
study $L^p$ spectral multipliers for operators without heat kernels. Our
approach provides a strong tool to verify the generalized Gaussian estimates
for a large class of operators. It is natural here to consider not only $L^p$
spaces but other functional spaces. This leads to another generalization of
Gaussian estimates (see Section~\ref{os} below).

Before we introduce all technical details needed to state our main results,
we would like to discuss Theorems~\ref{main2}~and~\ref{main22} below, which are
only specific consequences of these results, but provide a good
non-technical illustration of our approach. In \cite[Theorem~3.4.8 p.103]%
{Da} (see also \cite[Lemma~9]{Da1}) Davies shows that the Gaussian estimate
for the heat kernel extends to complex values of time. The surprisingly
simple proof of Theorem~\ref{main2}  yields a more precise version of 
\cite[Theorem~3.4.8 p.103]{Da}, and at the same time it provides an
alternative proof of real time off-diagonal Gaussian bounds obtained in \cite%
{Da, DP, Co1}. Recall that the heat semigroup $\exp (-t\Delta)$ generated by the
(non-negative) Laplace-Beltrami operator $\Delta$ on a complete Riemannian manifold $M$ is
self-adjoint on $L^{2}(M)$ with a smooth positive kernel $p_{t}(x,y)$, $t>0$%
, $x,y\in M$, called the heat kernel on $M$; it extends to a complex time
semigroup $\exp (-z\Delta)$, $z\in \mathbb{C}_{+}$, with a smooth kernel $%
p_{z}(x,y)$, $z\in \mathbb{C}_{+}$, $x,y\in M$. Denote by $d$ the geodesic distance on $M$.

\begin{theorem}
\label{main2} Let $p_{z}$, $z\in \mathbb{C}_{+}$, be the heat kernel on a
complete Riemannian manifold $M$. Suppose that 
\begin{equation}
p_{t}(x,x)\leq K{ t}^{-D/2},\ \forall \,t>0,\,x\in M,  \label{fuf2}
\end{equation}%
for some $K$and $D>0$. Then 
\begin{equation}
|p_{z}(x,y)|\leq eK({\mbox{\Small{\rm Re}}}z)^{-D/2}\left( 1+{%
\mbox{\Small{\rm Re}}}\frac{{d}^{2}(x,y)}{4z}\right) ^{D/2}\exp \left( -{%
\mbox{\Small{\rm Re}}}\frac{{d}^{2}(x,y)}{4z}\right)   \label{aaa2}
\end{equation}%
for all $z\in \mathbb{C}_{+}$, $x,y\in M$.
\end{theorem}

For $z=t \in \mathbb{R}_+$, estimates (\ref{aaa2}) can still be improved. It
is possible to prove that 
\begin{equation}
0\le p_t(x,y)\le Ct^{-D/2}\left(1+\frac{{d}^2(x,y)}{4t}\right)^{(D-1)/2}
\exp\left(-\frac{{d}^2(x,y)}{4t}\right)  \label{siko}
\end{equation}
(see \cite{Si1}). Moreover, it is known that the additional term $\left(1+%
\frac{{d}^2(x,y)}{4t}\right)^{(D-1)/2}$ cannot be removed in general from (%
\ref{siko}). See \cite{Mo} for a counterexample. However, using the
Phragm\'en-Lindel\"of technique we obtain the following variation of Theorem~%
\ref{main2}.

\begin{theorem}
\label{main22} Let $p_{z}$, $z\in \mathbb{C}_{+}$, be the heat kernel on a
complete Riemannian manifold $M$. Suppose that 
\begin{equation}
|p_{z}(x,y)|\leq K|z|^{-D/2},\ \forall \,z\in \mathbb{C}_{+},\,x,y\in M,
\label{fuf22}
\end{equation}%
for some $K$and$\ D>0$. Then 
\begin{equation*}
|p_{z}(x,y)|\leq eK|z|^{-D/2}\exp \left( -{\mbox{\Small{\rm Re}}}\frac{{d}%
^{2}(x,y)}{4z}\right)
\end{equation*}%
for all $z\in \mathbb{C}_{+}$, $x,y\in M$.
\end{theorem}

Theorems~\ref{main2}~and~\ref{main22} are straightforward consequences of
Theorems~\ref{truemain2}~and~\ref{truemain22} below, and the well-known fact
that the Laplace-Beltrami operator on complete Riemannian manifolds
satisfies Davies-Gaffney estimates, see the remark after Theorem~\ref{negpot}
below or \cite{Da2, Gr}.

Theorem~\ref{main22} shows that one can remove the additional factor $%
\left(1+\frac{{d}^2(x,y)}{t}\right)^{D/2}$ in (\ref{siko}) if one is able to
replace estimates (\ref{fuf2}) by the stronger ones (\ref{fuf22}). This is
an example of a result which we can obtain using Phragm\'en-Lindel\"of
technique and which does not seem to follow from the techniques developed in 
\cite{Da, DP, Co1, Gr0, Gr1, Gr, Si1, Si2}.

\section{Theorems of Phragm\'en-Lindel\"of type}

\label{pl}

Let us start with stating the Phragm\'en-Lindel\"of theorem for sectors.

\begin{theorem}
\label{phli} Let $S$ be the open region in $\mathbb{C}$ bounded by two rays
meeting at an angle $\pi/\alpha$, for some $\alpha>1/2$. Suppose that $F$ is
analytic on $S$, continuous on $\bar{S}$, and satisfies $|F(z)| \le C
\exp(c|z|^\beta)$ for some $\beta \in [0, \alpha) $ and for all $z\in S$.
Then the condition $|F(z)| \le B$ on the two bounding rays implies $|F(z)|
\le B$ for all $z\in S$.
\end{theorem}

For the proof see \cite[Theorem~7.5, p.214, vol.II]{Ma} or \cite[Lemma~4.2,
p.108]{SW}. Propositions~\ref{tw0},~\ref{tw1}~and~\ref{tw2} are simple
consequences of Theorem~\ref{phli}.

\begin{proposition}
\label{tw0} Suppose that $F$ is an analytic function on $\mathbb{C}_+$.
Assume that, for given numbers $A, B, \gamma>0$, $a\ge 0$, 
\begin{equation}  \label{a1}
|F(z)|\le B, \quad \forall \, z \in \mathbb{C}_+ , \qquad \mbox{and}
\end{equation}
\begin{equation}  \label{a2}
|F(t)|\le Ae^{at} e^{-\frac{\gamma}{t}}, \quad \forall \, t \in \mathbb{R}_+.
\end{equation}
Then 
\begin{equation}  \label{bb2}
|F(z)| \le B\exp \left(-{\mbox{\Small{\rm Re}}} \frac{\gamma}{z}\right),
\quad \forall \, z \in \mathbb{C}_+ .
\end{equation}
\end{proposition}

\begin{proof}
Consider the function 
\begin{equation}  \label{gg1}
u(\zeta)=F\left(\frac{\gamma}{\zeta}\right),
\end{equation}
which is also defined on $\mathbb{C}_+$. By (\ref{a1}), 
\begin{equation*}
|u(\zeta)e^\zeta|\le B \exp|\zeta|, \quad \forall \, \zeta \in \mathbb{C}_+ .
\end{equation*}
Again by (\ref{a1}) we have, for any $\varepsilon > 0$, 
\begin{equation}  \label{22a}
\qquad \sup_{{\mbox{\Small{\rm Re}}}\zeta=\varepsilon}|u(\zeta)e^\zeta|\le
Be^\varepsilon.
\end{equation}
By (\ref{a2}), 
\begin{equation}  \label{22b}
\sup_{ \zeta\in [\varepsilon,\infty)}|u(\zeta)e^\zeta|\le A
e^{a\gamma/\varepsilon}.
\end{equation}
Hence, by Phragm\'en-Lindel\"of theorem with angle $\pi/2$ and $\beta=1$,
applied to 
\begin{equation*}
S^+_{\varepsilon}= \{z\in \mathbb{C} \colon \, {\mbox{\Small{\rm Re}}} z >
\varepsilon \quad \mbox{and} \quad {\mbox{\small{\rm Im}}} z > 0\}
\end{equation*}
and 
\begin{equation*}
S^-_{\varepsilon}= \{z\in \mathbb{C} \colon \, {\mbox{\Small{\rm Re}}} z >\varepsilon
\quad \mbox{and} \quad {\mbox{\small{\rm Im}}} z < 0\},
\end{equation*}
one obtains 
\begin{equation*}
\sup_{{\mbox{\Small{\rm Re}}} \zeta \ge \varepsilon}|u(\zeta)e^\zeta|\le
\max\{ A e^{a\gamma/\varepsilon},Be^\varepsilon \}, \quad \forall \,
\varepsilon > 0 .
\end{equation*}
Now by the Phragm\'en-Lindel\"of theorem with angle $\pi$ and $\beta=0$, 
\begin{equation}  \label{22c}
\sup_{{\mbox{\Small{\rm Re}}} \zeta \ge \varepsilon}|u(\zeta)e^\zeta|\le
Be^\varepsilon, \quad \forall \, \varepsilon > 0 .
\end{equation}
Letting $\varepsilon \to 0$ we obtain 
\begin{equation*}
\sup_{{\mbox{\Small{\rm Re}}} \zeta > 0}|u(\zeta)e^\zeta|\le B.
\end{equation*}
This proves (\ref{bb2}) by putting $\zeta =\frac{\gamma}{z}$.
\end{proof}

Note that the estimate (\ref{bb2}) does not depend on constants $A,a$ in (%
\ref{a2}). This simple observation is the heart of the matter in the present
paper.

\medskip The above proposition will be used to prove the equivalence between
the finite speed propagation property for the solution of the wave equation
and Davies-Gaffney estimates (see \S \ref{dg} below). However, in order to
study the Gaussian bounds for heat kernels, we shall need a more
sophisticated version of Proposition~\ref{tw0}.

Given $\gamma>0$, denote by $\mathcal{C}_\gamma$ the closed disk in $\mathbb{%
C}_+$ centered on the real axis, tangent to the imaginary axis, with radius $%
\gamma/2$, that is the region 
\begin{equation*}
\mathcal{C}_\gamma=\{z\in \mathbb{C}\setminus\{0\}: {\mbox{\Small{\rm Re}}} 
\frac{\gamma}{z} \ge {1}\}.
\end{equation*}

\begin{proposition}
\label{tw1} Let $F$ be an analytic function on $\mathbb{C}_+$. Assume that,
for given numbers $A, B, \gamma,\nu>0$, 
\begin{equation}  \label{a11}
|F(z)|\le A
\end{equation}
for all $z\in \mathbb{C}_+$; 
\begin{equation}  \label{a22}
|F(t)|\le A e^{-\frac{\gamma}{t}}
\end{equation}
for all $t \in \mathbb{R}_+$ such that $t\le {\gamma} $; 
\begin{equation}  \label{a33}
|F(z)|\le B \left(\frac{{\mbox{\Small{\rm Re}}} z}{4\gamma}\right)^{-\nu/2}
\end{equation}
for all $z\in \mathcal{C}_\gamma$. Then 
\begin{equation}  \label{b1}
|F(z)| \le {e} B\left(\frac{2\gamma}{|z|}\right)^{\nu }\exp\left(- {%
\mbox{\Small{\rm Re}}}\frac{\gamma}{z} \right)
\end{equation}
for all $z\in \mathcal{C}_\gamma$.
\end{proposition}

\begin{proof}
Consider again the function $u$ defined by (\ref{gg1}). It satisfies
condition~(\ref{22a})~and~(\ref{22b}) with $B=A$, $a=0$ and $\varepsilon ={1}
$. Hence by (\ref{22c}) 
\begin{equation}  \label{we11}
\sup_{{\mbox{\Small{\rm Re}}} \zeta \ge {1}}|u(\zeta)e^\zeta|\le {e} A.
\end{equation}
Consider now the function $v$ defined on $\mathbb{C}_+$ by the formula 
\begin{equation*}
v(\zeta)= (2\zeta)^{-\nu} u(\zeta)e^\zeta.
\end{equation*}
Note that $|v(\zeta)|\le 2^{-\nu}|u(\zeta)e^\zeta|$ for ${%
\mbox{\Small{\rm
Re}}} \zeta \ge {1} $ so by (\ref{we11}) $v$ is bounded on the set ${%
\mbox{\Small{\rm Re}}} \zeta \ge {1} $ . Now, by Phragm\'en-Lindel\"of
theorem with angle $\pi$ and $\beta=0$, 
\begin{eqnarray*}
\sup_{{\mbox{\Small{\rm Re}}} \zeta \ge {1}} |v(\zeta)|= \sup_{{%
\mbox{\Small{\rm Re}}} \zeta = {1}} |v(\zeta)|.
\end{eqnarray*}
Put $\zeta={1}+is$. By (\ref{a33}), 
\begin{eqnarray*}
\sup_{{\mbox{\Small{\rm Re}}} \zeta ={1}} |v(\zeta)|&=& \sup_{{%
\mbox{\Small{\rm Re}}} \zeta ={1}} | (2\zeta)^{-\nu} F\left(\frac{\gamma}{%
\zeta}\right) e^\zeta| \\
&\le& \sup_{{\mbox{\Small{\rm Re}}} \zeta ={1}} {e} B |2\zeta|^{-\nu} \left({%
\mbox{\Small{\rm Re}}} \frac{1}{4\zeta}\right)^{-\nu/2} \\
&\le& {e} B \sup_{s\in \mathbb{R}} \,\, ({1+s^2})^{-\nu/2} \left(\frac{1}{{%
1+s^2}}\right)^{-\nu/2} = {e} B.
\end{eqnarray*}
Hence 
\begin{eqnarray}  \label{pl1}
\sup_{{\mbox{\Small{\rm Re}}} \zeta \ge {1}} |v(\zeta)|\le {e} B.
\end{eqnarray}
Now we put $\zeta=\frac{\gamma}{z}$ in (\ref{pl1}) and we obtain 
\begin{equation*}
\Big|\frac{2\gamma}{z}\Big|^{-\nu} \big|F(z)\big| \exp \left({%
\mbox{\Small{\rm Re}}} \frac{\gamma}{z}\right) \le {e} B .
\end{equation*}
This proves (\ref{b1}) for all $z\in \mathcal{C}_\gamma$.
\end{proof}

Finally let us discuss one more version of Proposition~\ref{tw1}. We shall
need this modified version to prove Theorem~\ref{main22}. The main
difference between Propositions~\ref{tw1}~and~\ref{tw2} is that we multiply
the function $F$ by an analytic function $e^g$ satisfying a growth condition.

\begin{proposition}
\label{tw2} Let $g$ and $F$ be analytic functions on $\mathbb{C}_+$. Assume
that, for given numbers $C,c, \gamma>0$ and $0\le \beta <1$, 
\begin{equation}  \label{a10}
\left|\exp g(z)\right|\le Ce^{c|z|%
^{-\beta}}
\end{equation}
for all $z\in \mathcal{C}_\gamma$.

Next assume that $F$ satisfies conditions \eqref{a11} and
\eqref{a22} and that 
\begin{equation}  \label{ab33}
|F(z)|\le B \left(\frac{{%
\mbox{\Small{\rm Re}}} z}{4\gamma}\right)^{-\nu/2}\exp(-{\mbox{\Small{\rm Re}}} {g(z)})
\end{equation}
for some $\nu>0$ and all $z\in \mathcal{C}_\gamma$. Then 
\begin{equation*}
|F(z)|\le {e} B\left(\frac{2\gamma}{|z|}%
\right)^{\nu} \exp\left(-{\mbox{\Small{\rm Re}}} {g(z)}- {\mbox{\Small{\rm Re}}}\frac{\gamma}{z} \right)
\end{equation*}
for all $z\in \mathcal{C}_\gamma$.
\end{proposition}

\begin{proof}
Define functions $u_1$ and $v_1$ on $\mathbb{C}_+$ by the formulae 
\begin{equation*}
v_1(\zeta)=(2\zeta)^{-\nu}
u_1(\zeta)e^\zeta=(2\zeta)^{-\nu}u(\zeta)\exp\left( g\left(\frac{\gamma}{%
\zeta}\right) \right) e^\zeta =(2\zeta)^{-\nu}F\left(\frac{\gamma}{\zeta}%
\right) \exp\left( g\left(\frac{\gamma}{\zeta}\right) \right)e^\zeta.
\end{equation*}
By (\ref{we11}) and (\ref{a10}), 
\begin{equation*}
|v_1(\zeta)|\le 2^{-\nu}|u_1(\zeta)e^{\zeta}| =2^{-\nu} | u(\zeta) e^\zeta|
| \exp\left( g\left(\frac{\gamma}{\zeta}\right) \right)| \le 2^{-\nu}
eACe^{c\gamma^{-\beta}|\zeta|^\beta}
\end{equation*}
if ${\mbox{\Small{\rm Re}}} \zeta \ge1$. Now by Phragm\'en-Lindel\"of
theorem with angle $\pi$ 
\begin{equation*}
\sup_{{\mbox{\Small{\rm Re}}} \zeta \ge {1}} |v_1(\zeta)|= \sup_{{%
\mbox{\Small{\rm Re}}} \zeta = {1}} |v_1(\zeta)|.
\end{equation*}
By (\ref{ab33}), 
\begin{equation*}
\sup_{{\mbox{\Small{\rm Re}}} \zeta ={1}} |v_1(\zeta)|\le \sup_{{%
\mbox{\Small{\rm Re}}} \zeta ={1}} {e} B |2\zeta|^{-\nu} \left({%
\mbox{\Small{\rm Re}}} \frac{1}{4\zeta}\right)^{-\nu/2} ,
\end{equation*}
and the rest of the proof is as in Proposition ~\ref{tw1}.
\end{proof}

\section{Davies-Gaffney estimates}

Let $(M,d,\mu)$ be a metric measure space, that is $\mu$
is a Borel measure with respect to the topology defined by the metric $d$. Next let $B(x,r)=\{y\in M,\, {d}(x,y)< r\}$ be the open ball
with center $x\in M$ and radius $r>0$. For $1\le p\le+\infty$, we denote the
norm of a function $f\in L^p(M,d\mu)$ by $\|f\|_p$, by $\langle .,. \rangle$
the scalar product in $L^2(M,d\mu)$, and if $T$ is a bounded linear operator from $%
L^p(M,d\mu)$ to $L^q(M,d\mu)$, $1\le p,q\le+\infty$, we write $\|T\|_{p\to q} $ for
the  operator norm of $T$.

Suppose that, for every $z\in \mathbb{C}_+$, $\Psi(z)$ is a bounded linear operator
acting on $L^2(M,d\mu)$ and that $\Psi(z)$ is an analytic function of $z$.
Assume in addition that 
\begin{equation}  \label{DG1}
\|\Psi(z)\|_{2\to 2} \le 1, \quad \forall \, z \in \mathbb{C}_+ .
\end{equation}
For $U_1, U_ 2 \subset M$ open subsets of $M$, let $d(U_1,U_2)=\inf_{x\in
U_1, y\in U_2}{d}(x,y)$. We say that the family $\{\Psi(z) \colon \, z\in 
\mathbb{C}_+\}$ satisfies the Davies-Gaffney estimate if 
\begin{equation}  \label{DG2}
| \langle \Psi(t)f_1,f_2\rangle | \le \exp\left(-\frac{r^2}{4t}\right)
\|f_1\|_{2} \|f_2\|_{2}
\end{equation}
for all $t >0$, $U_i \subset M$, \ $f_i\in L^2(U_i,d\mu)$, $i=1,2$ and $%
r=d(U_1,U_2)$. Note that we only assume that (\ref{DG2}) holds for positive
real $t$.

A slightly different form of Davies-Gaffney estimate is mostly considered in the literature
(see for instance \cite{Da2} or \cite{Gr}): in our notation, it reads 
\begin{equation}  \label{DG2g}
| \langle \Psi(t)\chi_{U_1},\chi_{U_2}\rangle | \le \exp\left(-\frac{r^2}{4t}%
\right) \sqrt{\mu(U_1)\mu(U_2)}
\end{equation}
where $\chi_U$ denotes the characteristic function of the set $U$. Of
course, \eqref{DG2g} follows from \eqref{DG2} by taking $f_1=\chi_{U_1}$ and 
$f_2=\chi_{U_2}$. Conversely, assume \eqref{DG2g} and let $f_i=\sum_j c^i_j
\chi_ {A_j^i}$, where $A_j^i \subset U_i$. Then 
\begin{eqnarray*}
<\Psi(t) f_1,f_2> &\le& \sum_j \sum_\ell |c^1_j c^2_\ell|
(\mu(A^1_j)\mu(A^2_\ell))^{1/2}\exp\left(-\frac{d^2(A^1_j,A^2_\ell)}{4t}%
\right) \\
&\le& \sum_j \sum_\ell |c^1_j c^2_l|
(\mu(A^1_j)\mu(A^2_\ell))^{1/2}\exp\left(-\frac{d^2(U_1,U_2)}{4t}\right).
\end{eqnarray*}
By \eqref{DG1}, $<\Psi(z) f_1,f_2> \le \|f_1\|_2\|f_2\|_2$. Proposition \ref%
{tw0} then yields \eqref{DG2} for such $f_1,f_2$, and one concludes by
density.

\bigskip

One may wonder what is the justification of the  constant 4
in (\ref{DG2}); we shall see in Theorem \ref{fspro} below that in the case
where $\Psi(z)$ is a semigroup $e^{-zL}$, 4 is the good normalisation
between the operator $L$ and the distance $d$, namely it translates the fact that the
associated wave equation has propagation speed 1.

The other constants in \eqref{DG1} and \eqref{DG2} have been normalized to
one for simplicity, anyway then can be absorbed by multiplying accordingly
the family $\Psi(z)$.

\subsubsection{Examples}

Semigroups of operators generated by non-negative self-adjoint  operators always
satisfy (\ref{DG1}), and among them many examples of interest satisfy (\ref%
{DG2}). Recall that, if $L$ is a non-negative self-adjoint  operator on $%
L^2(M,d\mu) $, one can construct the spectral decomposition $E_L(\lambda)$ of
the operator $L$. For any bounded Borel function $m\colon [0, \infty) \to 
\mathbb{C}$, one then defines the operator $m(L)\colon \, L^2(M,d\mu) \to
L^2(M,d\mu)$ by the formula 
\begin{equation*}
m(L)=\int_0^{\infty}m(\lambda) \,\mathrm{d} E_L(\lambda).
\end{equation*}
Now, for $z\in\mathbb{C}_+$ and $m_z(\lambda) =\exp(-z\lambda)$, one sets $%
m_z(L)=\exp(-zL)$, $z\in \mathbb{C}_+$. By spectral theory the family $%
\Psi(z)=\{\exp(-zL)\colon \, z \in \mathbb{C}_+\}$, also called semigroup of
operators generated by $L$, satisfies condition~(\ref{DG1}).

As we already said, condition~(\ref{DG2}) holds for  all kinds of
self-adjoint, elliptic, second order like operators. Condition~(\ref{DG2})
is well-known to hold for Laplace-Beltrami operators on all complete
Riemannian manifolds. More precisely, Condition~(\ref%
{DG2g}) is proved for such operators  in \cite{Da2} and \cite{Gr}.
See also the remark after Theorem \ref{negpot}. In the more general setting
of Laplace type operators acting on vector bundles, condition~(\ref{DG2}) is
proved in \cite{Si2}. Another important class of semigroups satisfying
condition~(\ref{DG2}) are semigroups generated by Schr\"odinger operators
with real potential and magnetic field (see for example \cite{Si3}, as well
as Theorem \ref{negpot} and Section \ref{vb} below).

Note that self-adjointness  and non-negativity of $L$ are a way to ensure that
$\exp{(-zL)}$ is defined for $z\in\mathbb{C}_+$, and (\ref{DG1}), but these conditions
may hold for non-self-adjoint $L$, and they are sufficient by themselves to run the rest of our theory.

Estimates~(\ref{DG2}) also hold in the setting of  local Dirichlet
forms (see for example  \cite[Theorem~2.8]{HR}, and also \cite{Stu1}, \cite{Stu2}). In this case the metric
measure spaces under consideration are possibly not equipped with any
differential structure. However, the  semigroups associated with these Dirichlet forms  do satisfy in general  Davies-Gaffney estimates with respect to an intrinsic distance.

In the sequel, if $(M,d,\mu)$ is a metric measure space and $L$ a
non-negative self-adjoint  operator on $L^2(M,d\mu)$, we shall say abusively that $%
(M,d,\mu, L)$ satisfies the Davies-Gaffney condition if ~\textrm{(\ref{DG2})}
holds with $\Psi(t)=e^{-tL}$.

\subsection{Self-improving properties of Davies-Gaffney estimates}

It is convenient to establish two simple lemmas concerning Davies-Gaffney
estimates before discussing our main results. First we observe that given 
\textrm{(\ref{DG1})} it is enough to test \textrm{(\ref{DG2})} on balls
only. Then we observe that any additional multiplicative constant or even
additional exponential factor in \textrm{(\ref{DG2})} can be replaced by the
constant in \textrm{(\ref{DG1}). }

\begin{lemma}
\textrm{\label{ext} Suppose that $(M,d, \mu)$ is a separable metric space and
that the analytic family $\{\Psi(z) \colon \, z\in \mathbb{C}_+\}$ of
bounded operators on $L^2(M,d\mu)$ satisfies condition \eqref{DG1}, and
condition \eqref{DG2} restricted to all balls $U_i= B(x_i,r_i)$, $i=1,2$,
for all $x_1,x_2\in M$, $r_1,r_2>0$. Then it satisfies condition  \eqref{DG2}  for all open subsets $U_1,U_2$. }
\end{lemma}

\begin{proof}
\textrm{Let $U_1$ and $U_2$ be arbitrary open subsets of $M$; set $%
r=d(U_1,U_2)$. Let $f=\sum_{i=1}^kf_{i}$, where for all $1\le i \le k$, $%
f_{i} \in L^2( B(x_{i},r_{i}), d\mu)$, $B(x_{i},r_{i})\subset U_1$, and $%
f_{i_1}(x)f_{i_2}(x)=0$ for all $x\in M$, $1\le i_1 < i_2 \le k$. Similarly
let $g=\sum_{j=1}^\ell g_{j}$ where $g_{j}\in L^2( B(y_{j},s_{j}), d\mu)$, $%
B(y_{j},s_{j})\subset U_2$ for all $1\le j\le \ell$, and $%
g_{j_1}(x)g_{j_2}(x)=0$ for all $x\in M$, $1\le j_1 < j_2 \le \ell$. Note
that $d(B(x_{i},r_{i}),B(y_{j},s_{j})) \ge r.$ Now if condition ~(\ref{DG2})
holds for balls then 
\begin{eqnarray*}
| \langle \Psi(t)f,g\rangle |& =&| \langle \Psi(t) \sum_{i=1}^kf_{i},
\sum_{j=1}^\ell g_{j} \rangle | \\
&=&\sum_{i=1}^k\sum_{j=1}^\ell| \langle \Psi(t) f_{i}, g_{j} \rangle | \\
&\le&\sum_{i=1}^k\sum_{j=1}^\ell e^{-\frac{r^2}{4t}}\|f_i\|_2 \|g_j\|_2 \\
&\le&e^{-\frac{r^2}{4t}}\left(\sum_{i=1}^k\|f_i\|_2\right)\left(\sum_{j=1}^%
\ell \|g_j\|_2\right) \\
&\le&e^{-\frac{r^2}{4t}}\sqrt{k\ell } \left(\sum_{i=1}^k\|f_i\|_2^2%
\right)^{1/2}\left(\sum_{j=1}^\ell \|g_j\|_2^2\right)^{1/2} \\
&=& e^{-\frac{r^2}{4t}}\sqrt{k\ell }\|f\|_2 \|g\|_2.
\end{eqnarray*}
Now we assume that (\ref{DG1}) holds so if we put $F(z)=\langle
\Psi(z)f,g\rangle$ then Proposition~\ref{tw0} shows that the term $Ckl$ in
the above inequality can be replaced by 1. This means that (\ref{DG2}) holds
for $f$ and $g$. Now to finish the proof of the lemma, it is enough to note
that, since $M$ is separable, the space of all possible finite linear
combinations of functions $f$ such that ${\rm supp} f \subset B(x,r)
\subset U$ is dense in $L^2(U, d\mu)$. Moreover, if $f=\sum_{i=1}^kf_i$ and $%
f_i\in L^2(B(x_i,r_i), d\mu)$ for all $1\le i\le k$ then there exist functions $%
\tilde{f}_i\in L^2(B(x_i,r_i), d\mu)$ such that $f=\sum_{i=1}^k\tilde{f}_i$ and in
addition, for all $1\le i_1 < i_2 \le k$, $\tilde{f}_{i_1}(x)\tilde{f}%
_{i_2}(x)=0$ for all $x\in M$.}
\end{proof}

\begin{lemma}
\textrm{\label{ex1} Suppose that the family $\{\Psi(z) \colon \, z\in 
\mathbb{C}_+\}$ satisfies condition \eqref{DG1}. Assume in addition that,
for some $C\ge 1$ and some $a>0$, 
\begin{equation}  \label{fs2weak}
| \langle \Psi(t)f_1,f_2\rangle | \le Ce^{at}e^{-\frac{r^2}{4t}} \|f_1\|_{2}
\|f_2\|_{2}, \quad \forall t >0,
\end{equation}
whenever $f_i\in L^2(M,d\mu)$, ${\rm supp} f_i\subseteq B(x_i,r_i)$, $i =
1,2 $, and $r={d}(B(x_1, r_1), B(x_2,r_2))$. Then the family $\{\Psi(z)
\colon \, z\in \mathbb{C}_+\}$ satisfies condition \eqref{DG2}. }
\end{lemma}

\begin{proof}
\textrm{Lemma~\ref{ex1} is a straightforward consequence of Proposition~\ref%
{tw0} and Lemma~\ref{ext}. }
\end{proof}

\textrm{Let us give an application of Lemma~\ref{ex1} by giving yet another
example where Davies-Gaffney estimates hold, namely Schr\"odinger semigroups
with real potential. Suppose that $\Delta$ is the non-negative Laplace-Beltrami operator on a
Riemannian manifold $M$ with Riemannian measure $\mu$ and geodesic distance $d$, and consider the operator $\Delta+\mathcal{V}$ acting on $%
C_c^\infty(M)$, where $\mathcal{V} \in L^1_{\mbox{ \tiny{\rm
loc}}}(M,d\mu)$. If we assume that $\Delta+\mathcal{V} \ge 0$ then we can define
the Friedrichs extension of $\Delta+\mathcal{V}$, which with some abuse of
notation we also denote by $\Delta+\mathcal{V}$ (see for example \cite[%
Theorem~1.2.8]{Da}. }

\begin{theorem}
\textrm{\label{negpot} Suppose that $\Delta$ is the Laplace-Beltrami operator on
a Riemannian manifold $M$, that $\mathcal{V} \in L^1_{\mbox{ \tiny{\rm loc}}%
}(M, d\mu)$ and that $\Delta+\mathcal{V} \ge 0$ as a quadratic form. Then the semigroup 
$\{\Psi(z)=\exp(-z(\Delta+\mathcal{V}))\colon \, z\in \mathbb{C}_+\}$ satisfies
condition \eqref{DG2}. }
\end{theorem}

\begin{proof}
We start our proof with the additional assumption $\mathcal{V} \ge 0$%
. For $f\in L^2(M,d\mu)$, $t>0$, $x\in M$, we put ${f}_t(x)={f}(t,x)=\exp(-t(\Delta+\mathcal{V})){f}(x)$. Let $\kappa>0$, and a
function $\xi\in C^\infty(M)$, both to be chosen later, such that $|\nabla \xi|\le \kappa$, where $%
\nabla$ is the Riemannian gradient on $M$. Next, as in \cite{Da1, Gr, Si2}, we consider the integral
 \begin{equation*}
E(t)=\int_M|{f}(t,x)|^2e^{\xi(x)}\,d\mu(x).
\end{equation*}
Then 
\begin{eqnarray*}
\frac{E^{\prime}(t)}{2}&=&{\mbox{\Small{\rm Re}}} \int_M\partial_t{f}(t,x)%
\overline{{f}(t,x)}e^{\xi(x)}\,d\mu(x) =-{\mbox{\Small{\rm Re}}}\int_M\left((\Delta+\mathcal{V}){f}_t\right)\,%
\overline{{f}_t }e^\xi\,d\mu \\
&=&-{\mbox{\Small{\rm Re}}} \int_M \left(\nabla {f}_t\cdot \nabla({f}_t e^\xi) +|{f%
}_t|^2\mathcal{V} e^\xi\right) \,d\mu \\
&=&-{\mbox{\Small{\rm Re}}} \int_M \left(|\nabla {f}_t|^2+\nabla {f}_t\cdot{f%
}_t \nabla \xi +|{f}_t|^2\mathcal{V}\right)e^\xi \,d\mu \\
&\le& \int_M \left(-|\nabla {f}_t|^2+|\nabla {f}_t| |\nabla \xi ||{f}_t |
\right)e^\xi \,d\mu \\
&\le& \frac{1}{4}\int_M |{f}_t|^2|\nabla\xi|^2 e^\xi\,d\mu\le \frac{\kappa^2
E(t)}{4}
\end{eqnarray*}
(note that the non-negativity of $\mathcal{V}$ is used in the first inequality).
Hence $E(t)\le \exp(\kappa^2t/2)E(0)$. 

Consider now two disjoints open sets $U_1$ and $U_2$ in $M$.   Choose $\xi=\kappa d(.,U_1)$.  One has $|\nabla
\xi| \le \kappa$,  $\xi\equiv 0$ on $ U_1$, and, for any $g\in L^2_{\mbox{ \tiny{\rm loc}}%
}(M,d\mu)$, 
\begin{equation*}
 \int_{U_2}|{g}|^2 e^{\xi}\,d\mu\ge
e^{\kappa r} \int_{U_2}|{g}|^2\,d\mu,
\end{equation*}
where $r = d(U_1,U_2) $. 
Hence if supp~$f \subseteq U_1$ then, taking $g=f_t$, 
\begin{equation*}
\int_{U_2}|{f}_t|^2\,d\mu \le e^{-\kappa r}E(t)\le \exp\left(\frac{\kappa^2t}{2}-\kappa r\right)E(0)= \exp\left(\frac{\kappa^2t}{2}-\kappa r\right)%
\int_{U_1}|f|^2\,d\mu .
\end{equation*}
Choosing finally $\kappa=r/t$ we obtain 
\begin{equation*}  
\int_{U_2}|\Psi(t)f|^2\,d\mu\le \exp\left(-\frac{r^2}{2t}\right)
\int_{U_1}|f|^2\,d\mu, 
\end{equation*}
that is, for all $f\in L^2(U_1,d\mu)$,
\begin{equation*} 
\sup_{g\in L^2(U_2,d\mu),\ \|g\|_2=1}|\langle \Psi(t)f,g\rangle | ^2= \int_{U_2}|\Psi (t) f|^2\,d\mu\le \exp\left(-\frac{r^2}{2t}\right)
\|f\|_2^2, 
\end{equation*}
which yields  (\ref{DG2}). Next, we consider a potential $\mathcal{V}\in L^1_{%
\mbox{
\tiny{\rm loc}}}(M)$ such that $\Delta+\mathcal{V}\ge 0$. We put $\mathcal{V}%
_a(x)=\max\{\mathcal{V}(x),-a\}$ and $L_a=\Delta+\mathcal{V}_a$. When $a$ goes to 
$\infty$ then $L_a$ converges to $L=\Delta+\mathcal{V}$ in the strong resolvent
sense (see \cite[Theorem~VIII.3.3, p.454]{Ka} or \cite[Theorem S.16 p.373]%
{RS}). Hence by \cite[Theorem~VIII.20, p.286]{RS} or by \cite[%
Theorem~VIII.3.11, p.459 and Theorem~IX.2.16, p.504 ]{Ka}, $\exp(-tL_a)f$
converges to $\exp(-tL)f=\exp(-t(\Delta+\mathcal{V}))f$ for any $f \in L^2(M)$. Hence it is
enough to prove (\ref{DG2}) for a given $a\in \mathbb{R}$. Finally we note
that $\mathcal{V}_a+a\ge 0$, thus it follows from the first part of the
proof that 
\begin{equation*}
\exp(-t(\Delta+\mathcal{V}_a+a))=e^{-at} \exp(-t(\Delta+\mathcal{V}_a))
\end{equation*}
satisfies condition~(\ref{DG2}). But this implies that the semigroup $%
\exp(-t(\Delta+\mathcal{V}_a))$ satisfies condition~(\ref{fs2weak}) and Theorem~%
\ref{negpot} follows from Lemma~\ref{ex1}. 
\end{proof}

\textrm{\emph{Remark :} Note that the case $\mathcal{V} =0$ is allowed in
Theorem~\ref{negpot}, in other words it yields a proof of (\ref{DG2}) for
the Laplace-Beltrami operator on complete Riemannian manifolds. }

\subsection{\textrm{Finite speed propagation for the wave equation and
Davies-Gaffney estimates}}

\textrm{\label{dg} }

\textrm{As a next application of the Phragm\'en-Lindel\"of technique
developed in Section \ref{pl}, we show that, for self-adjoint operators,
Davies-Gaffney estimates are equivalent to finite speed propagation property
for the corresponding wave equation. This equivalence was proved, along the
same lines, in \cite{Si2}. The underlying idea is so basic to the
development of the technique in the present paper that we shall repeat this
proof. One can use this equivalence to obtain a very simple proof of the
finite speed propagation property for a broad class of self-adjoint
operators (see \cite{Si2}). We start with recalling the notion of finite
speed propagation property for the wave equation. }

\textrm{In this section, $(M,d,\mu)$ is again a metric measure space. We say
that a  non-negative self-adjoint  operator $L$ satisfies the finite speed
propagation property for  solutions of the corresponding wave equation if 
\begin{equation}  \label{fs11}
\langle \cos(t\sqrt L) f_1, f_2 \rangle = 0
\end{equation}
for all $0<t < r$, open sets $U_i \subset M$, \ $f_i\in L^2(U_i,d\mu)$, $%
i=1,2 $, where $r=d(U_1,U_2)$. }

\textrm{If $\cos(t\sqrt L)$ is an integral operator with kernel $K_t$, then (%
\ref{fs11}) simply means that ${\rm supp} K_{t} \subseteq \mathcal{D}_t$,
that is $K_t(x_1,x_2)=0$ for all $(x_1,x_2) \notin \mathcal{D}_t$, where 
\begin{equation*}
\mathcal{D}_t=\{ (x_1,x_2)\in M\times M: {d}(x_1,x_2) \le t \}.
\end{equation*}
}

\begin{theorem}
\textrm{\label{fspro} Let $L$ be a self-adjoint non-negative operator acting on $%
L^2(M,d\mu)$. Then the finite speed propagation property  \eqref{fs11} and
Davies-Gaffney estimates \eqref{DG2} are equivalent. }
\end{theorem}

\begin{proof}
\textrm{Assume the Davies-Gaffney estimates. Fix two open sets $U_1,
U_2\subset M$. Let $f_i\in L^2(U_i, d\mu)$ for $i=1,2$. Define a function $%
F\colon \mathbb{C}_+ \to \mathbb{C}$ by 
\begin{equation}  \label{ab}
F(z)= \langle \exp(-zL)f_1,f_2\rangle.
\end{equation}
Since $\exp(-zL)$ is contractive on $L^2(M,d\mu)$, $F$ is a bounded analytic
function on $\mathbb{C}_+$ and it satisfies (\ref{a1}) with 
\begin{equation*}
B=\|f_1\|_{2} \|f_2\|_{2}.
\end{equation*}
In virtue of (\ref{DG2}), $F$ satisfies (\ref{a2}) with 
\begin{equation*}
a=0, \quad A=\|f_1\|_{2} \|f_2\|_{2}, \quad \gamma=\frac{r^2}{4} \quad %
\mbox{and} \quad r={d}(U_1,U_2),
\end{equation*}
thus, by Proposition~\ref%
{tw0}, 
\begin{equation}\label{biel}
\left|F(z)\right|\le \|f_1\|_{2} \|f_2\|_{2} \exp \left(-r^2{\mbox{\Small{\rm Re}}} \frac{1}{4z}\right).
\end{equation}
Then write, for $s>0$,  the well-known  Hadamard transmutation formula
\begin{equation}  \label{fc}
<\exp(-s L)f_1,f_2>= \int_0^{\infty} < \cos(t\sqrt L)f_1,f_2> \frac{e^{-%
\frac{t^2}{4s}}}{\sqrt {\pi s}} \,d t.
\end{equation}
By the change of variable $t\to\sqrt t$ in integral (\ref{fc}) and changing $%
s$ to $1/(4s)$, we obtain 
\begin{equation*}
{s}^{-1/2}<\exp{\left(-\frac{L}{4s}\right)}f_1,f_2>= \int_0^{\infty} ({\pi t}%
)^{-1/2} <\cos(\sqrt t \sqrt L)f_1, f_2> e^{-st} \,dt,
\end{equation*}
and by analytic continuation 
\begin{equation}  \label{dok}
{\zeta}^{-1/2}<\exp{\left(-\frac{L}{4\zeta}\right)}f_1,f_2>= \int_0^{\infty}
({\pi t})^{-1/2} <\cos(\sqrt t \sqrt L)f_1, f_2> e^{-\zeta t} \,dt
\end{equation}
for $\zeta\in\mathbb{C}_+$. Consider the function $u(\zeta)={\zeta}%
^{-1/2}F\left(\frac{1}{4\zeta}\right) $. On the one hand, by (\ref{biel}), 
\begin{equation*}
\left|u(\zeta)\right|\le \|f_1\|_{2} \|f_2\|_{2} |\zeta|^{-1/2}\exp (-r^2{\mbox{\Small{\rm Re}}} {\zeta}).
\end{equation*}
 On the other hand, by (\ref{dok}), $u$ is
the Fourier-Laplace transform of the function $$v(t)=(\pi t)^{-1/2} <
\cos(\sqrt t \sqrt L)f_1, f_2>.$$}

\textrm{A suitable version of the Paley-Wiener theorem (see \cite%
[Theorem~7.4.3, p.193]{H}) shows that 
\begin{equation*}
\supp\; v \subseteq [r^2, \infty).
\end{equation*}
Thus $< \cos(\sqrt t \sqrt L)f_1, f_2>=0$ for all $0\le t \le r^2$ and (\ref%
{fs11}) is proved. }

\textrm{Conversely if (\ref{fs11}) holds, then by (\ref{fc}) 
\begin{eqnarray*}
|< \exp(-s L)f_1,f_2>|&\le& \int_0^{\infty} \left|< \cos(t\sqrt
L)f_1,f_2>\right| \,\frac{e^{-\frac{t^2}{4s}} }{\sqrt {\pi s}} \,d t \\
&=& \int_r^{\infty} \left|< \cos(t\sqrt L)f_1,f_2>\right|\, \frac{e^{-\frac{%
t^2}{4s}} }{\sqrt {\pi s}} \,d t \\
&\le& \|f_1\|_{2}\|f_2\|_{2}\int_r^{\infty} \frac{e^{-\frac{t^2}{4s}}}{%
\sqrt {\pi s}} \,d t \\
&\le& e^{-\frac{r^2}{4s}}\|f_1\|_{2}\|f_2\|_{2}.
\end{eqnarray*}
}
\end{proof}

\section{\textrm{F{}rom on-diagonal bounds to Gaussian bounds}}

\textrm{\label{DGtoGB} }

\textrm{In this section, $(M,d,\mu)$ is again a metric measure space.
Let $L$ be a non-negative self-adjoint operator acting on $L^2(M,d\mu)$; recall that the
semigroup of operators $\Psi(z)=\exp(-zL)$, $z\in \mathbb{C}_+$,  is
contractive on $L^2(M,d\mu)$, in other words it satisfies condition~(\ref{DG1}%
). }
\medskip

Our basic observation is the following : if $p_z$ is the kernel associated with
$\exp(-zL)$, the estimate  $$p_t(x,x)\leq K{\ t}^{-D/2}, \ \forall\, t>0$$
can be reformulated as
\begin{equation*}
\|\exp(-tL)\|_{1 \to \infty}\le Kt^{-D/2}, \ \forall\, t>0
\end{equation*}
which yields
\begin{equation*}
\|\exp(-zL)\|_{1\to\infty}\le K({\mbox{\Small{\rm Re}}} z)^{-D/2}, \
\forall\, z\in \mathbb{C}_+
\end{equation*}
that is, in terms of a bilinear estimate,
\begin{equation*}
|<\exp(-zL)f_1,f_2>|\le K({\mbox{\Small{\rm Re}}} z)^{-D/2}\|f_1\|_1\|f_2\|_1, \
\forall\, z\in \mathbb{C}_+,  f_1,f_2\in L^1(M,d\mu).
\end{equation*}
On the other hand, Davies-Gaffney says that
\begin{equation*} 
| \langle \exp(-tL) f_1,f_2\rangle | \le \exp\left(-\frac{r^2}{4t}\right)
\|f_1\|_{2} \|f_2\|_{2}
\end{equation*}
for all $t >0$,  \ $f_1,f_2\in L^2(M,d\mu)$, supported respectively in $U_1,U_2$, with  $%
r=d(U_1,U_2)$.
With these ingredients and the global $L^2$ bound on the complex half-plane
\begin{equation*}
|<\exp(-zL)f_1,f_2>|\le \|f_1\|_2\|f_2\|_2, \
\forall\, z\in \mathbb{C}_+,  f_1,f_2\in L^2(M,d\mu), 
\end{equation*}
the suitable Phragm\'en-Lindel\"of  type lemma yields
\begin{equation*}
|\langle \exp (-zL)f_{1},f_{2}\rangle | \le eK({\mbox{\Small{\rm Re}}}z)^{-D/2}\left( {\mbox{\Small{\rm Re}}}\frac{%
r^{2}}{4z}\right) ^{D/2}\exp \left( -{\mbox{\Small{\rm Re}}}\frac{r^{2}}{4z}\right)\|f_1\|_1\|f_2\|_1,
\end{equation*}
for $z$ in the relevant region of the half-plane and  all  \ $f_1,f_2\in L^1(M,d\mu)$, supported  in $U_1,U_2$, with  $%
r=d(U_1,U_2)$.
Taking for $U_1$, $U_2$  balls that shrink around $x$ and $y$,
one obtains, for ${\mbox{\Small{\rm Re}}}\frac{%
d^{2}(x,y)}{4z}$ large enough, the desired estimate
\begin{equation*}
|p_z(x,y) | \le eK({\mbox{\Small{\rm Re}}}z)^{-D/2}\left( {\mbox{\Small{\rm Re}}}\frac{%
d^{2}(x,y)}{4z}\right) ^{D/2}\exp \left( -{\mbox{\Small{\rm Re}}}\frac{d^{2}(x,y)}{4z}
\right).
\end{equation*}

Let us see this in more detail, and then in more general situations.

\subsection{\textrm{Polynomial decay}}

\textrm{\label{poly} }

\textrm{Assume 
\begin{equation*}
\|\exp(-tL)\|_{1 \to \infty}\le Kt^{-D/2}, \ \forall\, t>0.
\end{equation*}
It follows that 
\begin{equation}
\|\exp(-zL)\|_{1\to\infty}\le K({\mbox{\Small{\rm Re}}} z)^{-D/2}, \
\forall\, z\in \mathbb{C}_+.  \label{ct}
\end{equation}
Indeed, for $t>0, \,s\in\mathbb{R}$, 
\begin{eqnarray*}
\|\exp(-(t+is)L)\|_{1\to \infty}&\le& \|\exp(-tL/2)\|_{1 \to 2}
\|\exp(-isL)\|_{2 \to 2}\|\exp(-tL/2)\|_{2 \to \infty} \\
&=& \|\exp(-tL/2)\|^2_{1 \to 2} \|\exp(-isL)\|_{2 \to 2} \\
&=& \|\exp(-tL)\|_{1 \to \infty} \|\exp(-isL)\|_{2 \to 2} \\
&\le& K {\ t}^{-D/2},
\end{eqnarray*}
using the well-known equality $\|T^*T\|_{1\to \infty}=\|T\|_{1\to 2}^2$. }

\textrm{In particular, by \cite[Theorem~6, p.503]{DS}, $\exp(-zL)$ is an
integral operator for all $z\in \mathbb{C}_+$. This means there exists a
measurable kernel, which we denote by $p_z(x,y)$, such that 
\begin{equation}  \label{kernel1}
[\exp(-zL)f](x)=\int_M p_z(x,y)f(y)\,d\mu(y), \mbox{ for a.e. }x\in M.
\end{equation}
Before we start discussing Gaussian bounds, let us state another
straightforward consequence of \cite[Theorem~6, p.503]{DS}, which we are
going to use frequently in the sequel. } A linear operator $S$ is bounded from 
$L^1(M,d\mu)$ to
$L^\infty(M,d\mu)$ if and only if it is an integral operator with kernel $p(x,y)$ such that
$\mbox{\rm esssup}_{x\in M,y\in M} |p(x,y)|$ is finite, in which case
$$\mbox{\rm esssup}_{x\in M,y\in M} |p(x,y)|=\|S\|_{1\to \infty}.$$
More precisely, we have the following :

\begin{cl}
\textrm{\label{ker1} Let $U_1,U_2$ be open subsets of $M$. If $
p(x,y)$ is the kernel of a linear operator $S\colon \, L^1(M,d\mu) \to
L^\infty(M,d\mu)$, then 
\begin{equation*}
\mbox{\rm esssup}_{x\in U_1,y\in U_2} |p(x,y)|= \sup\left\{ |\langle Sf_1,f_2\rangle|
\colon \, \|f_1\|_{L^1(U_1,d\mu)} = \|f_2\|_{L^1(U_2,d\mu)}=1 \right\}.
\end{equation*}
 }
\end{cl}

To complete the last step of the argument we sketched above, namely to pass from estimates on arbitrarily small balls to pointwise estimates, we need to assume  space continuity of the kernel under consideration, which is the case in most concrete situations, but not in general.
\textit{We shall assume from now on that for every $z\in \mathbb{C}_+$ the
kernel $p_z$ is a continuous complex-valued function defined on $M\times M$.}\footnote{For an interesting discussion about continuity properties of a general heat kernel, see \cite{Grf}. On the other hand, Brian Davies told us about a folklore example of a decent Schr\"odinger operator on $\R$ whose on-diagonal values of the kernel are null on a countable dense subset.}
As a consequence, we can replace the essential suprema by suprema in the above expressions,
and also record the following :

\textrm{If $p_t(x,y)$ is the kernel of $\exp(-tL)$, a well-known argument using the semigroup property and the fact that $p_t(y,x)={\overline{p_t(x,y)}}$  (see the proof of  \eqref{cam} below)
shows
further that }
\begin{equation}
\|\exp(-tL)\|_{1 \to \infty}=\sup_{x,y\in M}|p_t(x,y)|=\sup_{x\in M}p_t(x,x).\label{pus}
\end{equation}

\textrm{We can now state the general version of Theorem \ref{main2}. An even more general version will be given in Corollary \ref{ccoC} below, at the expense of a slightly more complicated proof. }

\begin{theorem}
\textrm{\label{truemain2} Assume that $(M,d,\mu ,L)$ satisfies the
Davies-Gaffney condition \eqref{DG2}. If, for some }$K\ $and \textrm{$D>0$, 
\begin{equation}
p_t(x,x)\leq K{\ t}^{-D/2},\ \forall
\,t>0, \ x\in M, \label{tfuf2}
\end{equation}%
then 
\begin{equation}
|p_{z}(x,y)|\leq eK({\mbox{\Small{\rm Re}}}z)^{-D/2}\left( 1+{%
\mbox{\Small{\rm Re}}}\frac{{d}^{2}(x,y)}{4z}\right) ^{D/2}\exp \left( -{{%
\mbox{\Small{\rm Re}}}}\frac{{d}^{2}(x,y)}{4z}\right)   \label{ge2}
\end{equation}%
for all $z\in \mathbb{C}_{+}$, $x,y\in M$. }
\end{theorem}

\begin{proof}
\textrm{Fix $x,y\in M$, and for ${d}(x,y)>2s>0$ define a bounded analytic
function $F\colon \mathbb{C}_{+}\rightarrow \mathbb{C}$ as in (\ref{ab}) by
the formula 
\begin{equation*}
F(z)=\langle \exp (-zL)f_{1},f_{2}\rangle ,
\end{equation*}%
where $f_{1}\in L^{1}(B(x,s), d\mu)\cap L^{2}(B(x,s), d\mu)$, $f_{2}\in
L^{1}(B(y,s), d\mu)\cap L^{2}(B(y,s), d\mu)$ and $\Vert f_{1}\Vert _{1}=\Vert f_{2}\Vert
_{1}=1$. In virtue of Davies-Gaffney estimates (\ref{DG2}) and (\ref{DG1}), $%
F$ satisfies (\ref{a11}) and (\ref{a22}) with 
\begin{equation*}
\gamma =r^{2}/4,\mbox{ where }r=d(x,y)-2s,\quad \mbox{and}\quad A=\Vert
f_{1}\Vert _{2}\Vert f_{2}\Vert _{2}<\infty .
\end{equation*}%
Assumption (\ref{tfuf2}) yields, through \eqref{pus} and \eqref{ct}, 
\begin{equation*}
|F(z)|=|\langle \exp (-zL)f_{1},f_{2}\rangle |\leq K({{\mbox{\Small{\rm Re}}}%
z})^{-D/2},\ \forall \,z\in \mathbb{C}_{+},
\end{equation*}%
so that $F$ satisfies (\ref{a33}) with $\nu =D$ and $B=Kr^{-D}$. By
Proposition~\ref{tw1}, 
\begin{eqnarray*}
|F(z)|=|\langle \exp (-zL)f_{1},f_{2}\rangle | &\leq &{e}Kr^{-D}\left( \frac{%
r^{2}}{2|z|}\right) ^{D}\exp \left( -{\mbox{\Small{\rm Re}}}\frac{r^{2}}{4z}%
\right) \\
&=&eK\left( \frac{r^{2}}{4|z|^{2}}\right) ^{D/2}\exp \left( -{%
\mbox{\Small{\rm
Re}}}\frac{r^{2}}{4z}\right) \\
&=&eK({\mbox{\Small{\rm Re}}}z)^{-D/2}\left( {\mbox{\Small{\rm Re}}}\frac{%
r^{2}}{4z}\right) ^{D/2}\exp \left( -{\mbox{\Small{\rm Re}}}\frac{r^{2}}{4z}%
\right)
\end{eqnarray*}%
for all $z\in \mathcal{C}_{r^{2}/4}$. }

\textrm{Hence by Claim~\ref{ker1} 
\begin{eqnarray*}
|p_z(x,y)|&\le& \sup_{x^{\prime}\in B(x,s)}\sup_{y^{\prime}\in
B(y,s)}|p_z(x^{\prime},y^{\prime})| \\
&=& \sup\left\{\langle \exp(-zL)f_1,f_2\rangle|\colon \,
\|f_1\|_{L^1(B(x,s), d\mu)} = \|f_2\|_{L^1(B(y,s), d\mu)}=1 \right\} \\
&\le& eK({\mbox{\Small{\rm Re}}} z)^{-D/2}\left({\mbox{\Small{\rm Re}}} 
\frac{r^2}{4z}\right)^{D/2}\exp\left(- {\mbox{\Small{\rm Re}}}\frac{r^2}{4z}%
\right)
\end{eqnarray*}
for all $z\in \mathcal{C}_{r^2/4}$. }

\textrm{Letting $s$ go to $0$ we obtain (\ref{ge2}) for ${%
\mbox{\Small{\rm
Re}}} \frac{{d}^2(x,y)}{4z}\ge {1}$. Finally for ${\mbox{\Small{\rm Re}}} 
\frac{{d}^2(x,y)}{4z}< {1}$, (\ref{ge2}) is a straightforward consequence of
(\ref{tfuf2}). Indeed in that case the Gaussian correction term satisfies $%
\exp\left(- {\mbox{\Small{\rm Re}}} \frac{{d}^2(x,y)}{4z}\right) > e^{-1}$,
and the estimate (\ref{ge2}) follows from 
\begin{equation*}
|p_z(x,y)| \le K({\mbox{\Small{\rm Re}}} z)^{-D/2}
\end{equation*}
which in turn follows from (\ref{ct}). }
\end{proof}

\textrm{Now for the general version of Theorem \ref{main22}. Again, an even more general version will be given in Theorem \ref{grri}  below.}

\begin{theorem}
\textrm{\label{truemain22} Assume that $(M,d,\mu ,L)$ satisfies the
Davies-Gaffney condition \eqref{DG2}. Next suppose that 
\begin{equation}
|p_{z}(x,y)|\leq K|z|^{-D/2},\ \forall \,z\in \mathbb{C}_{+},\,x,y\in M,  \label{tfuf22}
\end{equation}%
for some $K\ and\ D>0$. Then 
\begin{equation}
|p_{z}(x,y)|\leq eK|z|^{-D/2}\exp \left( -{\mbox{\Small{\rm Re}}}\frac{{d}%
^{2}(x,y)}{4z}\right)  \label{taaa22}
\end{equation}%
for all $z\in \mathbb{C}_{+}$, $x,y\in M$. }
\end{theorem}

\begin{proof}
\textrm{One starts as in the proof of Theorem \ref{truemain2}. Then  condition (\ref{tfuf22}) yields
\begin{equation*}
F(z):=|\langle \exp(-zL)f_1,f_2\rangle|\le K | z|^{-D/2}.
\end{equation*}
Choosing $g$ so that $|\exp(g(z))|=|z|^{D/2}$ and taking $B=K$,$\gamma=r^2/4$, $%
\beta=0$, $\nu=0$ in Proposition~\ref{tw2}, one obtains 
\begin{equation*}
|\langle \exp(-zL)f_1,f_2\rangle|\le K |z|^{-D/2} \exp\left(-{%
\mbox{\Small{\rm Re}}} \frac{r^2}{4z}\right)
\end{equation*}
and the rest of the proof is as before. }
\end{proof}

\textrm{\emph{Remarks :} The fact that an on-diagonal estimate for the heat
kernel implies an off-diagonal estimate is of course not new. See for
example \cite[Theorem~3.2.7, p.89]{Da} for the real time estimate and \cite[%
Theorem~3.4.8, p.103]{Da} for the complex time estimate. Note however, that
the results obtained in \cite{Da} are less precise than (\ref{ge2}) because
they involve $4+\varepsilon$ instead of $4$ in the exponential factor. To
our knowledge, the estimates (\ref{ge2}) with $4$ as an exponential factor
are new for complex time. On the other hand, for real time and diffusion
semigroups, estimates (\ref{ge2}) were obtained in \cite{DP} (see also \cite%
{Co1}). }

\textrm{For $z =t \in \mathbb{R}_+$, the estimates (\ref{ge2}) can still be
improved. It is possible to prove that 
\begin{equation*}
|p_t(x,y)| \le Ct^{-D/2}\left(1+\frac{{d}^2(x,y)}{4t}\right)^{(D-1)/2}
\exp\left(-\frac{{d}^2(x,y)}{4t}\right)
\end{equation*}
(see \cite{Si1}), and this is sharp due to \cite{Mo}. It is an interesting
question why our results here (and results in \cite{DP} and \cite{Co1}) give
weaker estimates with $D/2$ instead of $(D-1)/2$. It is so because in our
proof we do not use the fact that the family of operators under
consideration is a semigroup generated by a self-adjoint operator; for more
on this, see the discussion in the remark at the end of Section~\ref{nosemi}. 
}

\textrm{\bigskip }

\textrm{Suppose now that the self-adjoint contractive semigroup $\exp(-tL)$ on $L^2(M,d\mu)$ is in addition uniformly bounded on $L^\infty(M,d\mu)$,  which includes the case of the heat semigroup  on a
complete Riemannian manifold, since it is submarkovian.  Suppose also that estimates (\ref{tfuf2}) hold and that 
\begin{equation}  \label{shr}
\|\exp\left(isL\right)\|_{1\to\infty}\le C|s|^{-D/2}, \quad \forall \, s\in 
\mathbb{R}\setminus \{0\}.
\end{equation}
Then the semigroup $\exp(-zL)$ satisfies condition~(\ref{tfuf22}), hence the
corresponding heat kernel satisfies estimates~(\ref{taaa22}). Indeed, by (%
\ref{shr}) for all $t>0, s\in\mathbb{R}$, 
\begin{equation*}
\|\exp(-(t+is)L)\|_{1\to \infty}\le \|\exp(-tL/2)\|_{\infty \to \infty}
\|\exp(-isL)\|_{1 \to \infty} \le C {\ |s|}^{-D/2}.
\end{equation*}
Together with (\ref{ct}) this yields 
\begin{equation*}
\|\exp(-(t+is)L)\|_{1\to \infty}\le C\min\{t^{-D/2},|s|^{-D/2}\} \le
C^{\prime}|t+is|^{-D/2}
\end{equation*}
for all $t>0, s\in\mathbb{R}$ (as a matter of fact, (\ref{tfuf22}) is
equivalent to the conjunction of (\ref{tfuf2}) and (\ref{shr})). This shows
that Gaussian bounds without an additional polynomial correction factor are
a necessary condition for (\ref{shr}) to hold. Let us observe that estimates
(\ref{shr}) play an essential role in studying Strichartz type estimates
(see for example \cite{KT}). }

\subsection{\textrm{The doubling case}}

\textrm{Let $(M,{d},\mu)$ be a metric measure space as above, and let $p_z$, 
$z \in \mathbb{C}_+$ be a continuous heat kernel corresponding to a non-negative
self-adjoint operator $L$ on $L^2(M,d\mu)$. }

\textrm{One says that $(M,{d},\mu)$ satisfies the doubling property if there
exists $C>0$ such that 
\begin{equation}  \label{dd}
\mu(B(x, 2r)) \le C \mu(B(x, r)), \ \ \forall r> 0, \,x \in M.
\end{equation}
If this is the case, there exist $C, \delta>0$ such that 
\begin{equation}
\frac{\mu(B(x, s))}{\mu(B(x, r))} \le C \left(\frac{s}{r}\right)^\delta , \ \ \forall s \ge r> 0, \, x
\in M.  \label{ddt}
\end{equation}
}

\textrm{In such a situation, the most natural on-diagonal estimates for heat
kernels are of the type 
\begin{equation*}
p_t(x,x) \le \frac{C}{\mu(B(x,\sqrt t))}, \ \forall\,t>0, \,x\in M
\end{equation*}
(see for instance \cite{Gr}).}

\textrm{We are going to consider estimates of a similar form, but where the
quantity ${\mu(B(x,\sqrt t))}$ will be replaced by a function $V$ of $x$ and 
$t$ that is not necessarily connected with the volume of balls. }

\textrm{We shall assume that ${V}\colon M\times\mathbb{R}_+ \to \mathbb{R}_+$
is non-decreasing in the second variable, that is ${V}(x,s) \le {V}(x,r)$
for all $x\in M$ and all $0< s \le r$, and that it satisfies the doubling
condition 
\begin{equation}  \label{d}
\frac{{V}(x,s)}{{V}(x,r)}\le K' \left(\frac{s}{r}\right)^\delta 
\end{equation}
for all $s\ge r>0$ and all $x\in M$, and some constants $\delta \ge 0$ and $%
K' \ge 1$. Finally we shall assume that $V(x,t)$ is a continuous function of $%
x$.
}

\textrm{We shall then consider the on-diagonal estimate 
\begin{equation}  \label{on-dd}
p_t(x,x)V(x,\sqrt{t})\le 1, \ \forall\,x\in M,\, t>0.
\end{equation}
}
One should compare the following result with \cite[Theorem~4]{Si2}, which yields a slightly more precise estimate for real time.

\begin{theorem}
\textrm{\label{main1} Assume that $(M,d,\mu, L)$ satisfies the
Davies-Gaffney condition \eqref{DG2}. Next assume that the corresponding
heat kernel $p_z$ is continuous and satisfies the on-diagonal estimate \eqref%
{on-dd} with $V$ satisfying the doubling condition \eqref{d}. Then 
\begin{eqnarray}  \label{aaa1}
|p_z(x,y)| \le \frac{e K'} {\sqrt{V\left(x,\frac{{d}(x,y)}{2}\right)V\left(y,%
\frac{{d}(x,y)}{2}\right)} } \left(\frac{ {d}^{2}(x,y)}{4|z|}\right)^{\delta}
\exp\left(- {\mbox{\Small{\rm Re}}}\frac{{d}^2(x,y)}{4z}\right)
\end{eqnarray}
for all   $z\in \mathbb{C}_+$, $x,y\in M$ such that ${%
\mbox{\Small{\rm
Re}}}\frac{{d}^2(x,y)}{4z}\ge 1$. For all  $z\in \mathbb{C}_+$, $x,y\in M$,
and in particular if ${\mbox{\Small{\rm Re}}}\frac{{d}^2(x,y)}{4z}< 1$, one
has 
\begin{eqnarray}  \label{xy}
|p_{z}(x,y)| \le \frac{1}{\sqrt{V(x,\sqrt{{\mbox{\Small{\rm Re}}} z}) {V}(y,%
\sqrt{{\mbox{\Small{\rm Re}}} z})} }.
\end{eqnarray}
}
\end{theorem}

\begin{proof}
\textrm{For all $z\in\mathbb{C}_+$, $x,y\in M$, one has 
\begin{equation}\label{cam}
|p_{z}(x,y)| \le \sqrt{p_{{\mbox{\Small{\rm Re}}} z}(x,x) p_{{%
\mbox{\Small{\rm Re}}} z}(y,y)}.
\end{equation}
Indeed, 
\begin{eqnarray*}
|p_{z}(x,y)|&=&\left|\int_M p_{z/2}(x,u)p_{z/2}(u,y)\,d\mu(u)\right| \\
&\le&\left(\int_M |p_{z/2}(x,u)|^2\,d\mu(u)\right)^{1/2}\left(\int_M
|p_{z/2}(u,y)|^2\,d\mu(u)\right)^{1/2} \\
&=&\left(\int_M p_{z/2}(x,u)p_{\bar{z}/2}(u,x)\,d\mu(u)\right)^{1/2}\left(%
\int_M p_{\bar{z}/2}(y,u)p_{z/2}(u,y)\,d\mu(u)\right)^{1/2} \\
&\le&\sqrt{p_{{\mbox{\Small{\rm Re}}} z}(x,x) p_{{\mbox{\Small{\rm Re}}}
z}(y,y)}.
\end{eqnarray*}
}

\textrm{In the second equality above, we have used the fact that, since $L$ is self-adjoint, $\overline{p_z(x,y)}=p_{\overline z}(y,x)$. Together with \eqref{on-dd}, this yields \eqref{xy}. }

\medskip

\textrm{For a function ${W} \colon M \to \mathbb{C}$, we denote by $\mathbf{m%
}_W$ the operator of multiplication by ${W}$, that is 
\begin{equation*}
(\mathbf{m}_W f)(x)={W}(x)f(x),
\end{equation*}
and if ${W}\colon M \times \mathbb{R}_+ \to \mathbb{C}$ then for $r\in 
\mathbb{R}_+$ we set 
\begin{equation*}
(\mathbf{m}_{{W}(\, \cdot \,,r)}f)(x)={W}(x,r)f(x).
\end{equation*}
}

\textrm{Let us now set $W(x,r)=\sqrt{V(x,r)}$. Similarly as in the proof of
Theorem~\ref{truemain2}, fix $x,y\in M$ and for ${d}(x,y)>2s>0$ set $r={d}%
(x,y)-2s$. Then define a bounded analytic function $F\colon \mathbb{C}%
_{+}\rightarrow \mathbb{C}$ by the formula 
\begin{equation}
F(z)=\langle \exp (-zL)\mathbf{m}_{W\left( \,\cdot \,,\frac{r}{2}\right)
}f_{1},\mathbf{m}_{W\left( \,\cdot \,,\frac{r}{2}\right) }f_{2}\rangle ,
\label{abnew}
\end{equation}%
where $$f_{1}\in L^{1}(B(x,s),d\mu)\cap L^{2}(B(x,s),V(\,\cdot \,,r/2)\,d\mu ),
f_{2}\in L^{1}(B(y,s),d\mu)\cap L^{2}(B(y,s),V(\,\cdot \,,r/2)\,d\mu ),$$ and $%
\Vert f_{1}\Vert _{1}=\Vert f_{2}\Vert _{1}=1$. In virtue of (\ref{DG2}) and
(\ref{DG1}), $F$ satisfies (\ref{a11}) and (\ref{a22}) with 
\begin{equation*}
A=\Vert \mathbf{m}_{W\left( \,\cdot \,,\frac{r}{2}\right) }f_{1}\Vert
_{2}\Vert \mathbf{m}_{W\left( \,\cdot \,,\frac{r}{2}\right) }f_{2}\Vert
_{2}<\infty \mbox{ and }\gamma =r^{2}/4.
\end{equation*}%
Note that, for  $z\in \mathcal{C}_{\gamma }$,
\begin{equation}
\gamma \geq \left( {\mbox{\Small{\rm Re}}}\frac{1}{z}\right) ^{-1}\geq {%
\mbox{\Small{\rm Re}}}z,
\label{real}
\end{equation}
hence $r/2\ge \sqrt{\mbox{\Small{\rm Re}}z}$.
Now by (\ref{on-dd}), \eqref{real} and (\ref{d}), \vspace{4pt} 
\begin{eqnarray*}
|F(z)| &=&\left\vert \int_{M}\int_{M}p_{z}(x^{\prime },y^{\prime }){W}\left(
y^{\prime },\frac{r}{2}\right) f_{1}(y^{\prime }){W}\left( x^{\prime },\frac{%
r}{2}\right) f_{2}(x^{\prime })\,d\mu (y^{\prime })\,d\mu (x^{\prime
})\right\vert \\
&\leq &\sup_{x^{\prime },y^{\prime }\in M}{W}\left( x^{\prime },\frac{r}{2}%
\right) |p_{z}(x^{\prime },y^{\prime })|{W}\left( y^{\prime },\frac{r}{2}%
\right) \\
&\leq &\sup_{x^{\prime },y^{\prime }\in M}\sqrt{\frac{V\left( x^{\prime },%
\frac{r}{2}\right) V\left( y^{\prime },\frac{r}{2}\right) }{V(x^{\prime },%
\sqrt{{\mbox{\Small{\rm Re}}}z})V(y^{\prime },\sqrt{{\mbox{\Small{\rm Re}}}z}%
)}}\leq K'\left( \frac{r}{2\sqrt{{\mbox{\Small{\rm Re}}}z}}\right) ^{\delta
}=K'2^{-\delta }\left( \frac{{\mbox{\Small{\rm Re}}}z}{r^{2}}\right)
^{-\delta /2}
\end{eqnarray*}%
so that $F$ satisfies (\ref{a33}) with $B=K'2^{-\delta }$ and $\nu =\delta $.
By Proposition~\ref{tw1}, 
\begin{equation*}
|F(z)|\leq {e}K'2^{-\delta }\left( \frac{r^{2}}{2|z|}\right) ^{\delta }\exp
\left( -{\mbox{\Small{\rm Re}}}\frac{r^{2}}{4z}\right)
\end{equation*}%
for all $z\in \mathcal{C}_{r^{2}/4}$. Note that $L^{1}(B(x,s),d\mu)\cap
L^{2}(B(x,s),V(\,\cdot \,,r/2)\,d\mu )$ is dense in $L^{1}(B(x,s),d\mu)$, so by
Claim~\ref{ker1}, 
\begin{eqnarray*}
{V}^{1/2}\left( x,\frac{r}{2}\right) |p_{z}(x,y)|{V}^{1/2}\left( y,\frac{r}{2%
}\right) &\leq &\sup_{x^{\prime }\in B(x,s)}\sup_{y^{\prime }\in B(y,s)}{V}%
^{1/2}\left( x^{\prime },\frac{r}{2}\right) |p_{z}(x^{\prime },y^{\prime })|{%
V}^{1/2}\left( y^{\prime },\frac{r}{2}\right) \\
&=&\sup \left\{ |F(z)|\colon \,\Vert f_{1}\Vert _{L^{1}(B(x,s),d\mu)}=\Vert
f_{2}\Vert _{L^{1}(B(y,s),d\mu)}=1\right\} \\
&\leq &eK'\left( \frac{r^{2}}{4|z|}\right) ^{\delta }\exp \left( -{%
\mbox{\Small{\rm
Re}}}\frac{r^{2}}{4z}\right)
\end{eqnarray*}%
for all $z\in \mathcal{C}_{r^{2}/4}$. Letting $s\rightarrow 0$, we obtain
the estimate \eqref{aaa1} for $z\in \mathcal{C}_{{d}^{2}(x,y)/4}$. }
\end{proof}

Note that taking
$$V(x,r)=K^{-1}r^D,\ r>0,x\in M,$$
one sees that Theorem \ref{truemain2} is a particular case of Theorem \ref{main1}.

\textrm{The estimate in the following corollary is less precise than the one
in Theorem~\ref{main1}, but its algebraic form is convenient for
calculations and it is enough for most applications; also, it can be
compared with the case $m=2$ of the estimates in \cite[Proposition~4.1]{CCO}%
. The improvement with respect to \cite{CCO} is that the constant inside the
exponential is optimal, at the expense of a necessary polynomial correction
factor. To state the result, we put $\theta=\arg z$ for all $z\in \mathbb{C}%
_+$, so that $\cos\theta=\frac{{\mbox{\Small{\rm Re}}} z}{|z|}$. }

\begin{coro}
\textrm{\label{ccoC} Under the assumptions of Theorem $\rm{\ref{main1}}$, 
\begin{equation}  \label{cco}
|p_z(x,y)| \le \frac{eK'\left(1+{\mbox{\Small{\rm Re}}} \frac{{d}^2(x,y)}{4z}\right)^{\delta}} {\sqrt{{V}\left(x,\sqrt{\frac{|z|}{\cos{\theta}}}\right) {V%
}\left(y,\sqrt{\frac{|z|}{\cos{\theta}}}\right)} } \exp\left(- {%
\mbox{\Small{\rm Re}}} \frac{{d}^2(x,y)}{4z}\right) \frac{1}{(\cos
\theta)^{\delta}}
\end{equation}
for all  $z\in \mathbb{C}_+$, $x,y\in M$. }
\end{coro}

\begin{proof}
\textrm{Note that ${\mbox{\Small{\rm Re}}} z^{-1} =|z|^{-1} \cos \theta $, hence 
\begin{equation*}
\frac{ d^2(x,y)}{4|z|}\cos \theta = {\mbox{\Small{\rm Re}}} \frac{%
{d}^2(x,y)}{4z}, \ z\in\mathbb{C}_+.
\end{equation*}
Moreover, if $z\in \mathcal{C}_{{d}^2(x,y)/4}$, 
\begin{equation*}
\frac{{d}(x,y)}{2}\ge ({\mbox{\Small{\rm Re}}} z^{-1})^{-1/2}=(|z|/\cos{%
\theta})^{1/2}
\end{equation*}
hence
\begin{equation*}
{V}\left(x,\frac{{d}(x,y)}{2}\right) \ge {V}\left(x,\sqrt{\frac{|z|}{\cos{%
\theta}}}\right).
\end{equation*}
Therefore
\begin{eqnarray*}
&&\frac{1}{\sqrt{{V}\left(x,\frac{{d}(x,y)}{2}\right){V}\left(y,\frac{{d}%
(x,y)}{2}\right)}} \left(\frac{ {d}^{2}(x,y)}{4|z|}\right)^{\delta}\exp\left(-{%
\mbox{\Small{\rm Re}}} \frac{{d}^2(x,y)}{4z}\right) \\
&\le& \frac{1}{\sqrt{{V}\left(x,\sqrt{\frac{|z|}{\cos{\theta}}}\right) {V}%
\left(y,\sqrt{\frac{|z|}{\cos{\theta}}}\right)} } \left({%
\mbox{\Small{\rm
Re}}} \frac{{d}^2(x,y)}{4z}\right)^{\delta} \exp\left(- {%
\mbox{\Small{\rm
Re}}} \frac{{d}^2(x,y)}{4z}\right) \frac{1}{(\cos \theta)^{\delta}}
\end{eqnarray*}
and for $z\in \mathcal{C}_{{d}^2(x,y)/4}$, (\ref{cco}) follows from %
\eqref{aaa1}. Finally we note that by (\ref{d}) 
\begin{equation*}
{V}\left(x,\sqrt{\frac{|z|}{\cos{\theta}}}\right)\le K'\,V(x,\sqrt{ |z|\cos{%
\theta}})(\cos{\theta})^{-\delta}= K'\,V(x,\sqrt{{\mbox{\Small{\rm Re}}} z}%
)(\cos{\theta})^{-\delta}
\end{equation*}
so that, for $z\notin \mathcal{C}_{{d}^2(x,y)/4}$, since $\exp\left(- {%
\mbox{\Small{\rm Re}}} \frac{{d}^2(x,y)}{4z}\right)\ge e^{-1}$, (\ref{cco})
is a straightforward consequence of (\ref{xy}). }
\end{proof}

\bigskip

It is certainly an interesting feature of  Corollary \ref{ccoC} that it yields estimates valid for time ranging in the whole  right half-plane, and that it does not require $V$ to be tied to the volume.
Let us however observe the following particular case of our result, for real time and estimates involving the volume growth function. It also follows from \cite[Proposition~ 5.2]{Gr0},
but our proof is  more direct, as it does not go through a Faber-Krahn type inequality.

\begin{coro}
\label{grigo} Let $p_{t}$, $t>0$, be the heat kernel on a
complete Riemannian manifold $M$, with Riemannian measure $\mu$ and geodesic distance $d$. 
Let $V(x,r)$ denote
$\mu(B(x,r))$, for $r>0$, $x\in M$. Assume that $M$ satisfies the doubling property, more precisely 
 let $K',\delta>0$ be such that $\eqref{d}$ is satisfied.  
 Suppose that 
\begin{equation}
p_{t}(x,x)\leq \frac{K}{V(x,\sqrt{t})},\ \forall \,t>0,\,x\in M,  \label{gfuf2}
\end{equation}%
for some $K>0$. Then 
\begin{equation*}  
p_t(x,y) \le \frac{eK'K\left(1+ \frac{d^2(x,y)}{4t}\right)^{\delta}} {\sqrt{V\left(x,\sqrt{t}\right) V%
\left(y,\sqrt{t}\right) }} \exp\left(- \frac{{d}^2(x,y)}{4t}\right) \end{equation*}
for all  $t>0$, $x,y\in M$. 
\end{coro}

Let us now consider the case where the heat kernel satisfies upper and lower estimates 
of the type \begin{equation}
 p_t(x,y)\simeq\frac {1}{\mu\left(B(x,t^{1/\beta})\right)}\exp\left(-\left(\frac{d^\beta(x,y)}{t}\right)^{\frac{1}{\beta-1}}\right),
\ \forall\,t>0, \,x,y\in M.\label{douwintro}
\end{equation}
This may happen when $(M,d,\mu, L)$ is a fractal space, endowed with a natural metric,  measure
and Laplacian,   for all values of $\beta$ between $2$ and $\delta+1$, where $\delta$ is  the exponent  in the doubling condition \eqref{d}; see for instance \cite{BSF}.
In such situations, usually, $\beta>2$. 
Let us now choose $V(x,t)=\mu\left(B(x,t^{2/\beta}\right)$, which is obviously a doubling function.
F{}rom \eqref{douwintro}, $p_t$ satisfies  \eqref{on-dd}, but it cannot satisfy
\begin{equation*}  
p_t(x,y) \le \frac{C} {V\left(x,\sqrt{t}\right) } \exp\left(- c\frac{{d}^2(x,y)}{t}\right)=
\frac{C} {\mu\left(B(x,t^{1/\beta})\right) } \exp\left(- \frac{{d}^2(x,y)}{4t}\right),
 \end{equation*}
since this  is   not compatible with the lower bound in \eqref{douwintro}. In view of Corollary \ref{ccoC}, the only possible conclusion is that  such a space $(M,d,\mu, L)$ does not satisfy Davies-Gaffney estimates, nor, according to \eqref{fspro}, the finite speed propagation property for the wave equation. We owe this remark to Alexander Teplyaev. There is no contradiction with the fact that  local Dirichlet forms do give rise to Davies-Gaffney estimates with respect to an intrinsic distance : in the case of fractals, this distance  degenerates, see the discussion in \cite[Section 3.2]{HR}.

\textrm{\bigskip }

\textrm{Finally let us discuss one more version of pointwise Gaussian
estimates. Here we do not need to consider any kind of doubling property.
The following result has some similarity with Corollary~5.5 of \cite{Gr}, in
the sense that, in assumption \eqref{ggr} below, $x,y$ do not range in the
whole space $M$, but only in two fixed regions $U_1,U_2$. However,
restricting our assumption to two fixed points $x,y$ as in \cite{Gr} seems
to raise technical difficulties that we are not going to face here. }

\begin{theorem}
\textrm{\label{grri} Assume that $(M,d,\mu, L)$ satisfies the Davies-Gaffney
condition \eqref{DG2}. Let $U_1$ and $U_2$ be open subsets of $M$. Suppose
that 
\begin{equation}  \label{ggr}
|p_z(x,y)|\le \exp(-{\mbox{\Small{\rm Re}}}
g(z)), \quad \forall \, z\in \mathbb{C}_+,\,  x\in U_1,\, y\in U_2,
\end{equation}
where  $g$ is analytic on $\mathbb{C}_+$ and 
satisfies the growth condition \eqref{a10} with $\gamma=\frac{r^2}{4}$ and $r=d(U_1,U_2)$. Then 
\begin{equation*}
|p_z(x,y)|\le \exp\left(1-{\mbox{\Small{\rm Re}}%
} g(z) - {\mbox{\Small{\rm Re}}}\frac{r^2}{4z}\right), \quad \forall \,
z\in\mathbb{C}_+,\, x\in U_1,\, y\in U_2.
\end{equation*}
}
\end{theorem}

\textrm{\emph{Remarks:}} 

- It may look surprising that the growth constraint on $g$ depends on $U_1,U_2$. This may be understood as follows : suppose a factor
$ \exp\left( - {\mbox{\Small{\rm Re}}}\frac{r^2}{4z}\right)$ is already present in estimate
 \eqref{ggr}  (which corresponds to $\beta=1$, a situation hopefully forbidden by \eqref{a10}) ; then one can certainly not multiply again the estimate by this factor!
 
 - Theorem \ref{grri} is a generalization of Theorem \ref{truemain22}, as one can see by taking
 $g(z)=\frac{D}{2}\log z-\log K$. 
 
 -In principle, one could use Theorem  \ref{grri} to add a Gaussian factor
 to estimates
 of the form $$|p_z(x,y)|\le \frac{1}{|V(x,y,z)|}, \quad \forall \, z\in \mathbb{C}_+,\,  x, y\in M,$$
 where $V$ is analytic in $z$ with a certain uniformity in $x,y$. We will not
 pursue this direction because of the lack of relevant examples.
 
 -Note that our result allows to some extent rapid growth at zero. In particular, it might be interesting to investigate the connection with \cite[Theorem~3.1]{BSC}.

\begin{proof}
\textrm{Once again we follow the idea of the proof of Theorem~\ref{truemain2}
and define a bounded analytic function $F\colon \mathbb{C}_+ \to \mathbb{C}$
by the formula 
\begin{equation*}
F(z)= \langle \exp(-zL)f_1,f_2\rangle,
\end{equation*}
where $f_i \in L^2(U_i, d\mu) \cap L^1(U_i, d\mu)$ and $\|f_1\|_1=\|f\|_2=1$. In virtue
of (\ref{DG2}) and (\ref{DG1}), $F$ satisfies (\ref{a11}) and (\ref{a22})
with 
\begin{equation*}
r=d(U_1,U_2), \quad \gamma=r^2/4 \quad \mbox{and}\quad A= \| f_1\|_{2} \|
f_2 \|_{2}< \infty.
\end{equation*}
Next, by assumption (\ref{ggr}), 
\begin{equation*}
|F(z)| \le |\exp(-g(z))|\|f_1\|_{L^1(U_1, d\mu)}\|f_2\|_{L^1(U_2, d\mu)} =\exp(-{%
\mbox{\Small{\rm Re}}} g(z)), \quad \forall \, z\in \mathcal{C}_{r^2/4},
\end{equation*}
that is, $F$ satisfies (\ref{ab33}) with $B=1$, $\nu=0$. 
By Proposition~\ref{tw2}, 
\begin{equation*}
|F(z)|\le \exp\left(-{\mbox{\Small{\rm Re}}} g(z)+1- {%
\mbox{\Small{\rm Re}}} \frac{{r}^2}{4z}\right), \quad \forall \, z\in 
\mathcal{C}_{r^2/4}.
\end{equation*}
Finally by Claim~\ref{ker1} 
\begin{eqnarray*}
\sup_{x\in U_1}\sup_{x\in U_2}|p_z(x,y)| &=& \sup\left\{\langle
\exp(-zL)f_1,f_2\rangle\colon \, \|f_1\|_{L^1(U_1, d\mu)} = \|f_2\|_{L^1(U_2, d\mu)}=1
\right\} \\
&\le & \exp\left(-{\mbox{\Small{\rm Re}}} g(z)+1 - {%
\mbox{\Small{\rm
Re}}} \frac{{r}^2}{4z}\right)
\end{eqnarray*}
for all $z\in \mathcal{C}_{r^2/4}$. The estimate for $z\notin \mathcal{C}%
_{r^2/4}$ follows directly from \eqref{ggr}. }
\end{proof}

\subsection{\textrm{Operators acting on vector bundles}}

\textrm{\label{vb} }

\textrm{Our approach works not only for operators acting on functions but
can also be applied to operators acting on vector bundles. To discuss the
vector bundle version of our results we need some additional notation. }

\textrm{Let $(M,d,\mu)$ be a metric measure space and suppose that $TM$ is a
continuous vector bundle with  base $M$, fibers $T_xM\simeq\mathbb{C}^n$ and with
continuous (with respect to $x$) scalar product $(\, \cdot \,, \, \cdot
\,)_x $ on $T_xM$. For $f(x)\in T_xM$ we put $|f(x)|_x^2=(f(x), f(x))_x$. To
simplify the notation, we will write $(\, \cdot \,, \, \cdot \,)$ and $|\,
\cdot \,|$ instead of $(\, \cdot \,, \, \cdot \,)_x$ and $|\, \cdot \,|_x$.
Now for sections $f$ and $g$ of $TM$ we put 
\begin{equation*}
\|f\|_{L^p(M, d\mu ; TM)}^p=\int_M|f(x)|^p\,d\mu(x) \qquad \mbox{and} \quad \langle
f,g \rangle= \int_M (f(x),g(x))\,d\mu(x).
\end{equation*}
}

\textrm{Now let us describe the notion of \emph{integral operators} for
vector bundles. For any point $(x,y) \in M^2$ we consider the space $%
T_y^*\otimes T_x$. The space $T_y^*\otimes T_x$ is canonically isomorphic to 
${\mbox { Hom}}\,(T_y,T_x)$, the space of all linear homeomorphisms from $%
T_y$ to $T_x$. Denote again by $|\, \cdot \,|$ the operator norm on $%
T_y^*\otimes T_x$ with respect to the norms $|\, \cdot \,|_x$ and $|\, \cdot
\,|_y$. }

\textrm{By $(T^*\otimes T) M^2$ we denote the continuous bundle with 
base space equal to $M^2$ and with  fiber over the point $(x,y)$ equal to 
$T_y^*\otimes T_x$. If there is a section $\vec{p}$ of $(T^*\otimes T) M^2$
such that $|\vec{p}|$ is a locally integrable function on $(M^2,\mu\times\mu)$
and $\vec{S}f_1$ is a section of $TM$ such that 
\begin{equation*}
\langle \vec{S}f_1,f_2\rangle = \int_{M} (\vec{S}f_1(x) ,f_2(x)) \,d \mu(x) = \int_{M}
(\vec{p}(x,y) \, f_1(y),f_2(x)) \,d \mu(y) \,d \mu(x)
\end{equation*}
for all sections $f_1$ and $f_2$ in~$C_c(TM)$, then we say that $\vec{S}$ is an 
\emph{integral operator} on sections of $TM$ with kernel $\vec{p}$.   As in the scalar case,
$\vec{S}$ is a bounded linear operator from $L^1(M, d\mu ; TM)$ to~$%
L^\infty(M, d\mu ; TM)$ if and only if $\vec{S}$ is an integral operator with kernel $\vec{p}$ such that
$\mbox{\rm esssup}_{x,y\in M} |\vec{p}(x,y)|$ is finite, and 
\begin{equation*}
\mbox{\rm esssup}_{x,y\in M} |\vec{p}(x,y)|= \| \vec{S}\|_{1\to \infty}.
\end{equation*}}
One also has the following vector-valued version of Claim \ref{ker1} :
 \begin{cl}
\textrm{\label{ker2} Let  $U_1,U_2$ be open subsets of $M$. If $
\vec{p}(x,y)$ is the kernel of a linear operator $\vec{S}\colon \, L^1(M,d\mu ; TM) \to
L^\infty(M,d\mu ; TM)$, then 
\begin{equation*}
\mbox{\rm esssup}_{x\in U_1,y\in U_2} |\vec{p}(x,y)|= \sup\left\{ |\langle \vec{S}f_1,f_2\rangle|
\colon \, \|f_1\|_{L^1(U_1,d\mu ; TM)} = \|f_2\|_{L^1(U_2,d\mu ; TM)}=1 \right\}.
\end{equation*}
}
\end{cl}

\textrm{Let us describe an example of Hodge type operator which generates a
semigroup satisfying conditions~(\ref{DG1}) and (\ref{DG2}) and acts on
vector bundles of $k$-forms on Riemannian manifolds. Suppose that $M$ is a
complete $n$-dimensional Riemannian manifold and $\mu$ is an absolutely continuous measure
with a smooth density not equal to zero at any point of $M$. By $%
\Lambda^kT^*M$, $k=0,...,n$, we denote the bundle of $k$-forms on $M$. For fixed $\beta,
\beta^{\prime}\in L^2(\Lambda^1T^*M)$ and $\gamma \in L^2(\Lambda^kT^*M)$,
we define the operator $\vec{L}=\vec{L}_{\beta,\beta^{\prime},\gamma}$ acting on $%
L^2(\Lambda^kT^*M$) by the formula 
\begin{equation}  \label{hod1}
\langle \vec{L} \omega, \omega\rangle = \int_M \left(|d_k\omega +\omega \wedge \beta|^2+
|d_{n-k}*\omega +*\omega \wedge \beta^{\prime}|^2+ |*\omega \wedge \gamma|^2
\right)\,d\mu,
\end{equation}
where $\omega$ is a smooth compactly supported $k$-form and $*$ is the Hodge
star operator. With some abuse of notation we also denote by $\vec{L}$ its
Friedrichs extension. Note that for example the Hodge-Laplace operator and
Schr\"odinger operators with real potentials and electromagnetic fields can
be defined by (\ref{hod1}). The following theorem was proved in \cite{Si2}. }

\begin{theorem}
\textrm{\label{fsp} The self-adjoint semigroup $\{\exp(-z\vec{L})\colon \, z\in \mathbb{C}_+\}$
generated by the operator $\vec{L}$ defined by \eqref{hod1} acts on $%
L^2(\Lambda^kT^*M)$ and satisfies \eqref{DG1} and \eqref{DG2}. }
\end{theorem}

\textrm{Theorems~\ref{truemain2},~\ref%
{truemain22},~\ref{main1},~\ref{grri} and Corollary~\ref{ccoC} can be
extended to the above setting of  operators acting on vector bundles.   For
example we can state Theorem~\ref{main1} in this setting in the
following way.  Again, compare with \cite[Corollary~9]{Si2}, which yields a slightly better estimate for real time
and the Hodge-Laplace operator (and more generally operators defined by (\ref{hod1})), but does not treat complex time. In what follows, $ \mbox{\rm Tr}$ denotes the trace of an endomorphism on a finite dimensional linear space. }

\begin{theorem}\label{doudou}
\textrm{  Let $(M,d,\mu)$ be a metric measure space endowed with a vector bundle $TM$  as above. Let $\vec{L}$ be a non-negative self-adjoint operator acting on $L^2(M,d\mu ; TM)$. Assume that $\vec{p}_z(x,y)$ is a continuous function of $x,y\in M$. Denote by $\vec{p}_z$, $z\in \mathbb{C}_+$, the kernel of
$\exp(-z\vec{L})$.
Let  
 $V \colon \mathbb{R}_+\times M \to \mathbb{R}_+$ be
a continuous function 
satisfying condition \eqref{d}.  Assume  that
\begin{equation}
V\left(x,\sqrt{t}\right)\mbox{\rm Tr }\vec{p}_t(x,x) \le 1, \quad \forall\, t>0,\, x\in M.\label{ooo47p}
\end{equation}
Then 
\begin{equation}\label{oo47}
|\vec{p}_z(x,y)| \le \frac{e K'} 
{\sqrt{V\left(x,\frac{{d}(x,y)}{2}\right)V\left(y,%
\frac{{d}(x,y)}{2}\right)}}  
\left(\frac{ {d}^{2}(x,y)}{4|z|}\right)^{\delta}
\exp\left(- {\mbox{\Small{\rm Re}}}\frac{{d}^2(x,y)}{4z}\right)
\end{equation}
for all $z\in \mathbb{C}_+$, $x,y\in M$ such that ${%
\mbox{\Small{\rm
Re}}}\frac{{d}^2(x,y)}{4z}\ge 1$.}
\end{theorem} 

{\bf Remark :}  Of course, one can transform the above estimate in a similar way as in Corollary
\ref{ccoC}.

\begin{proof} Note that the self-adjointness of $\exp(-t\vec{L})$ implies $p_{t}(y,x)=[p_{t}(x,y)]^*$. 
Denote by $|.|_{HS}$ the Hilbert-Schmidt norm of a linear operator.
One can write
\begin{eqnarray*}
\mbox{\rm Tr }\vec{p}_t(x,x)&=&\mbox{\rm Tr} \int p_{t/2}(x,y)p_{t/2}(y,x)\,d\mu(y)\\&=&\mbox{\rm Tr} \int p_{t/2}(x,y)[p_{t/2}(x,y)]^*\,d\mu(y)\\
 &=& \int  \mbox{\rm Tr}\left( p_{t/2}(x,y)[p_{t/2}(x,y)]^*\right)\,d\mu(y),
\end{eqnarray*}
thus
\begin{equation}\label{tra}
\mbox{\rm Tr }\vec{p}_t(x,x)=\int |\vec{p}_{t/2}(x,y)|_{HS}^2 \,d\mu(y).
\end{equation}

On the other hand, 
\begin{eqnarray*}
|\vec{p}_t (x,y)|^2  &\le&\int |\vec{p}_{t/2}(x,z)|^2 \,d\mu(z)  \int |\vec{p}_{t/2}(z,y)|^2 \,d\mu(z)\\ &\le&  \int |\vec{p}_{t/2}(x,z)|_{HS}^2 \,d\mu(z)  \int |\vec{p}_{t/2}(z,y)|_{HS}^2 \,d\mu(z), 
\end{eqnarray*}
since
$|.|\le |.|_{HS}$;
hence, using \eqref{tra} and \eqref{ooo47p},
\begin{equation}
|\vec{p}_t (x,y)|^2  \le \frac{1}{V(x,t^{1/2})V(y,t^{1/2})},\label{trap}
\end{equation}
that is
\begin{equation*}
\|\mathbf{m}_{{W}(\, \cdot \,, \sqrt{t})}\exp(-t\vec{L}) 
\mathbf{m}_{ {W}(\, \cdot \,, \sqrt{t})}\|_{1 \to \infty}\le 1,
\end{equation*}
where
 $W(x,r)=\sqrt{V(x,r)}$, $x\in M$, $r>0$.
 This estimate can be extended to complex times.
Indeed, one also has
$$\int |\vec{p}_{t}(x,z)|^2 \,d\mu(z)\le \int |\vec{p}_{t}(x,z)|_{HS}^2 \,d\mu(z)\le \frac{1}{ V(x,t^{1/2})},$$
that is
\begin{equation}
 \| \mathbf{m}_{W(\, \cdot \,, \sqrt{t})}\exp(-t\vec{L}) \|_{2 \to \infty}=\| \exp(-t\vec{L})\mathbf{m}_{W(\, \cdot \,, \sqrt{t})} \|_{1 \to 2}\le
1, \quad \forall t>0.\label{ooo47}
\end{equation}

Using  the contractivity of $\exp(-is\vec{L})$, $s\in\mathbb{R}$%
, on $L^2$, one has 
\begin{eqnarray*}
&&\|\mathbf{m}_{{W}(\, \cdot \,, \sqrt{{\mbox{\Small{\rm Re}}} z})}\exp(-z\vec{L}) 
\mathbf{m}_{ {W}(\, \cdot \,, \sqrt{{\mbox{\Small{\rm Re}}} z})}\|_{1 \to \infty}
\\
&\le&\|\mathbf{m}_{{W}(\, \cdot \,, \sqrt{{\mbox{\Small{\rm Re}}} z}%
)}\exp(-(z/2)\vec{L})\|_{2 \to \infty} \|\exp(-(z/2)\vec{L})\mathbf{m}_{{W}(\, \cdot \,, 
\sqrt{{\mbox{\Small{\rm Re}}} z})}\|_{1 \to 2} \\
&\le& \|\mathbf{m}_{{W}(\, \cdot \,, \sqrt{{\mbox{\Small{\rm Re}}} z}%
)}\exp(-({\mbox{\Small{\rm Re}}} z/2)\vec{L})\|_{2 \to \infty} \|\exp(-({%
\mbox{\Small{\rm Re}}} z/2)\vec{L})\mathbf{m}_{{W}(\, \cdot \,, \sqrt{{%
\mbox{\Small{\rm Re}}} z})}\|_{1 \to 2} \\
&=& \|\exp(-({\mbox{\Small{\rm Re}}} z/2)\vec{L})\mathbf{m}_{{W}(\, \cdot
\,, \sqrt{{\mbox{\Small{\rm Re}}} z})}\|_{1 \to 2}^2.
\end{eqnarray*}
Together with  \eqref{ooo47} and the identity $\|T^*T\|_{1\to
\infty}=\|T\|^2_{1\to 2}$, this yields
\begin{equation}\label{o47}
\|\mathbf{m}_{{W}(\, \cdot \,, \sqrt{{\mbox{\Small{\rm Re}}} z})}\exp(-z\vec{L}) 
\mathbf{m}_{ {W}(\, \cdot \,, \sqrt{{\mbox{\Small{\rm Re}}} z})}\|_{1 \to \infty}\le 1,\ \forall\,z\in\mathbb{C}_+.
\end{equation}

\textrm{Similarly as in (\ref{abnew}), fix $x,y\in M$ and for ${d}(x,y)>2s>0$ set $r={d}%
(x,y)-2s$.  Consider the function $F$ defined
by the formula 
\begin{equation*}
F(z)=\langle \exp(-z\vec{L}) \mathbf{m}_{W\left(\, \cdot \,,\frac{r}{2}%
\right)}\omega_2,
\mathbf{m}_{W\left(\, \cdot \,,\frac{r}{2}\right)}\omega_1\rangle,
\end{equation*}
with $\omega_1 \in L^1(B(x,s),d\mu ; TM)\cap L^2(B(x,s),V(\, \cdot \,,r/2)\,d\mu ; TM)$,  $\omega_2 \in L^1(B(y,s),d\mu ; TM)\cap   L^2(B(y,s),V(\, \cdot \,,r/2)\,d\mu ; TM)$,
 and $\|\omega_1\|_{1}= \|\omega_2\|_{1}=1$. }
\textrm{In virtue of assumption~(\ref{DG1}) and Davies-Gaffney estimates (%
\ref{DG2}), $F$ satisfies (\ref{a11}) and (\ref{a22}) with $\gamma=r^2/4$
and 
\begin{equation*}
A= \| \mathbf{m}_{W\left(\, \cdot \,, \frac{r}{2}\right)}\omega_1\|_{2} \|%
\mathbf{m}_{W\left(\, \cdot \,, \frac{r}{2}\right)}\omega_2
\|_{2}=\|\omega_1\|_{L^2(U_1,V(\, \cdot \,,r/2)\,d\mu)} \|\omega_2
\|_{L^2(U_2,V(\, \cdot \,,r/2)\,d\mu)}.
\end{equation*}
Now if $z\in \mathcal{C}_{r^2/4}$, then $\sqrt{{%
\mbox{\Small{\rm Re}}} z}\le r/2$ by (\ref{real}). Using the assumptions on $\omega_1,\omega_2$, $%
W$ as well as (\ref{o47}), we obtain 
\begin{eqnarray*}
|F(z)| &=& |\langle \mathbf{m}_{W\left(\, \cdot \,, \frac{r}{2}%
\right)}\exp(-z\vec{L})\mathbf{m}_{W\left(\, \cdot \,, \frac{r}{2}\right)}\omega_2 ,
\omega_1 \rangle | \\
&\le& \sup_{x,y\in M} \frac{W\left(x, \frac{r}{2}\right) W\left(y, \frac{%
r}{2}\right) }{{W}(x,\sqrt{{\mbox{\Small{\rm Re}}} z}) {W}(y,\sqrt{{%
\mbox{\Small{\rm Re}}} z}) }|\langle \mathbf{m}_{ {W}(\, \cdot \,,\sqrt{{%
\mbox{\Small{\rm Re}}} z})}\exp(-z\vec{L})\mathbf{m}_{ {W}(\, \cdot \,,\sqrt{{%
\mbox{\Small{\rm Re}}} z})}\omega_2 ,
\omega_1 \rangle |  \\
&\le& \sup_{x,y\in M} \frac{W\left(x, \frac{r}{2}\right) W\left(y, \frac{%
r}{2}\right) }{{W}(x,\sqrt{{\mbox{\Small{\rm Re}}} z}) {W}(y,\sqrt{{%
\mbox{\Small{\rm Re}}} z}) } \| \mathbf{m}_{ {W}(\, \cdot \,,\sqrt{{%
\mbox{\Small{\rm Re}}} z})} \exp(-z\vec{L})\mathbf{m}_{ {W}(\, \cdot \,,\sqrt{{%
\mbox{\Small{\rm Re}}} z}) }\|_{1 \to \infty} \\
&\le& K'\left(\frac{r^2}{{4{\mbox{\Small{\rm Re}}} z}}%
\right)^{\delta/2}
\end{eqnarray*}
for all $z\in \mathcal{C}_{r^2/4}$. Thus $F$ satisfies (\ref{a33}) with $%
B=K' 2^{-\delta}$ and $\nu=\delta$. By
Proposition~\ref{tw1}, 
\begin{equation*}
|F(z)| \le {e}K'  \left(\frac{r^2}{4|z|}%
\right)^{\delta }\exp\left(- {\mbox{\Small{\rm Re}}}\frac{r^2}{4z},%
\right)
\end{equation*}
for all $z\in \mathcal{C}_{r^2/4}$. One finishes the proof as in Theorem \ref{main1}, using Claim \ref{ker2}.}
\end{proof}

\subsection{\textrm{Gaussian estimates for the gradient of the heat kernel}}

\textrm{The technique which we developed above can be applied to obtain
Gaussian bounds for gradient of the heat kernels. The following result is
motivated by some considerations in \cite{ACDH}, Section 1.4. In particular,
it is proved in \cite{ACDH} that under the assumptions below, for  $a=1/2$, the Riesz
transform is bounded on $L^p(M,d\mu)$ for $2<p<+\infty$. In the langage of 
\cite{ACDH}, we will show now that, under $(FK)$, conditions (1.7), (1.8) and $(G)$ are all equivalent, which was left open there. A similar result was obtained independently in
\cite{Du}, by a different method, relying directly on the finite speed propagation property for the wave equation. }

\begin{theorem}
\textrm{\label{gl} Let $M$ be a complete Riemannian manifold such that the
Riemannian measure $\mu$ satisfies the doubling condition \eqref{dd} and let 
$\Delta$ be the Laplace-Beltrami operator, $p_z$ the corresponding heat
kernel, $\nabla$ the Riemannian gradient on $M$. Suppose next that 
\begin{equation}  \label{on-d}
p_t(x,x) \le \frac{C}{\mu(B(x,\sqrt t))}, \ \forall\,t>0, \,x\in M,
\end{equation}
and that 
\begin{equation}  \label{lemgl2n}
\sup_{x,y\in
M}|\nabla p_t(x,y)|\mu(B(y,\sqrt{t})) \le Ct^{-a}
\end{equation}
for  some $a>0$ and all $t\in \mathbb{R}_+$. Then 
\begin{equation}  \label{rt1}
|\nabla p_t(x,y)| \le \frac {C}{t^a\mu(B(y, \sqrt t))} \left(1+\frac{{d}%
^2(x,y)}{4t}\right)^{3\delta+2a } \exp\left({-\frac{{d}^2(x,y)}{4t}}\right)
\end{equation}
for all $t>0$, $x,y\in M$. }
\end{theorem}

\textrm{In the proof of Theorem~\ref{gl} we shall need the following 
consequence of Corollary \ref{ccoC} and assumption. }

\begin{lemma}\label{lemgln}
\textrm{  Assume \eqref{dd},  \eqref{on-d}, and \eqref{lemgl2n}. 
Then 
\begin{equation}  \label{lemgl1n}
\|\nabla\exp(-z\Delta)\mathbf{m}_{V(\, \cdot \,,r)}\|_{1\to \infty}=\sup_{x,y\in
M}|\nabla p_z(x,y)|V(y,r) \le C({\mbox{\Small{\rm Re}}} z)^{-a}\left(\frac{|z|}{{\mbox{\Small{\rm Re}}} z}%
\right)^{2\delta}\left(\frac{r^2}{{\mbox{\Small{\rm Re}}} z}\right)^{\delta/2}
\end{equation}
for all $r>0$, $z\in \mathcal{C}_{r^2/4}$, where $\delta>0$ is the exponent in %
\eqref{ddt}. }
\end{lemma}

\begin{proof}
\textrm{
An immediate reformulation of \eqref{cco} is
\begin{equation*} 
|p_z(x,y)| \le \frac{eK'\left(1+{\mbox{\Small{\rm Re}}} \frac{{d}^2(x,y)}{4z}\right)^{\delta}} {\sqrt{{V}\left(x,\left({\mbox{\Small{\rm Re}}} \frac{1}{z}\right)^{-1/2}\right) {V%
}\left(y,\left({\mbox{\Small{\rm Re}}} \frac{1}{z}\right)^{-1/2}\right)} } \exp\left(- {%
\mbox{\Small{\rm Re}}} \frac{{d}^2(x,y)}{4z}\right)   \left(\frac{|z|}{{\mbox{\Small{\rm Re}}} z}%
\right)^{\delta},
\end{equation*}
which yields, for $0<c<1/4$,
\begin{equation*}  
|p_z(x,y)| \le \frac{C} {\sqrt{{V}\left(x,\left({\mbox{\Small{\rm Re}}} \frac{1}{z}\right)^{-1/2}\right) {V%
}\left(y,\left({\mbox{\Small{\rm Re}}} \frac{1}{z}\right)^{-1/2}\right)} } \exp\left(-c\, {%
\mbox{\Small{\rm Re}}} \frac{{d}^2(x,y)}{z}\right)   \left(\frac{|z|}{{\mbox{\Small{\rm Re}}} z}%
\right)^{\delta},
\end{equation*}
and, by doubling,
\begin{equation*}\label{nabo} 
|p_z(x,y)| \le \frac{C'} { {V%
}\left(y,\left({\mbox{\Small{\rm Re}}} \frac{1}{z}\right)^{-1/2}\right) } \exp\left(- c'{%
\mbox{\Small{\rm Re}}} \frac{{d}^2(x,y)}{z}\right)   \left(\frac{|z|}{{\mbox{\Small{\rm Re}}} z}%
\right)^{\delta}
\end{equation*}
Take now $z=t+is\in \mathcal{C%
}_{r^2/4}$, where $t,s\in \mathbb{R}$. Note that  $(t/2)+is \in \mathcal{C}_{r^2}$. Write
\begin{equation*}
\nabla p_z(x,y)=\int_M \nabla p_{t/2}(x,u) p_{(t/2)+is}(u,y)\,d\mu(u),
\end{equation*}
hence, using \eqref{nabo}, \eqref{lemgl2n}, and doubling,
\begin{equation*}
|\nabla p_z(x,y)|\le \frac{C}{t^{a} V\left(y,\left({\mbox{\Small{\rm Re}}} \frac{1}{z}\right)^{-1/2}\right) }\left(\frac{|z|}{{\mbox{\Small{\rm Re}}} z}\right)^{\delta}
\int_M \frac{1}{V(u,\sqrt{t})}   \exp\left(- c'{\mbox{\Small{\rm Re}}} \frac{{d}^2(u,y)}{z}\right)\,d\mu(u).
\end{equation*}
Let us estimate
\begin{equation*}\label{esint}
I=\int_M \frac{1}{V(u,\sqrt{t})}   \exp\left(- c'{\mbox{\Small{\rm Re}}} \frac{{d}^2(u,y)}{z}\right)\,d\mu(u).
\end{equation*}
Since $ \left({\mbox{\Small{\rm Re}}} \frac{1}{z}\right)^{-1/2}\ge {%
(\mbox{\Small{\rm Re}}}z)^{1/2}=\sqrt{t}$, one has
\begin{equation*}
I\le K'\left(\frac{({\mbox{\Small{\rm Re}}}\frac{1}{z})^{-1/2}}{\sqrt{{\mbox{\Small{\rm Re}}}z}}\right)^{\delta/2}II=\left( \frac{|z|}{{\mbox{\Small{\rm Re}}}z}\right)^\delta II,\end{equation*}
where
\begin{equation*}
II=\int_M \frac{1}{V\left(u,\left({\mbox{\Small{\rm Re}}}\frac{1}{z}\right)^{-1/2}\right)}   \exp\left(- c'{\mbox{\Small{\rm Re}}} \frac{{d}^2(u,y)}{z}\right)\,d\mu(u)
\end{equation*}
is easily seen to be uniformly bounded in $y\in M$, $z\in \mathbb{C}_+$ by doubling.
Thus
\begin{equation*}
|\nabla p_z(x,y)|\le \frac{C}{t^{a} V\left(y,\left({\mbox{\Small{\rm Re}}} \frac{1}{z}\right)^{-1/2}\right) }\left(\frac{|z|}{{\mbox{\Small{\rm Re}}} z}\right)^{2\delta}.
\end{equation*}
Now, for $z \in \mathcal{C}_{r^2}$,  $r\ge \left({\mbox{\Small{\rm Re}}} \frac{1}{z}\right)^{-1/2}\ge {%
(\mbox{\Small{\rm Re}}}z)^{1/2}$, hence
\begin{equation*}
\frac{V(y,r)}{V\left(y,\left({\mbox{\Small{\rm Re}}} \frac{1}{z}\right)^{-1/2}\right) }\le \frac{V(y,r)}{V\left(y,{%
(\mbox{\Small{\rm Re}}}z)^{1/2}\right) }\le K'\left(\frac{r^2}{{\mbox{\Small{\rm Re}}} z}\right)^{\delta/2}, 
\end{equation*}
therefore
\begin{equation*}  
|\nabla p_z(x,y)|V(y,r) \le C\left({\mbox{\Small{\rm Re}}} z\right)^{-a}\left(\frac{r^2}{{\mbox{\Small{\rm Re}}} z}\right)^{\delta/2}   \left(\frac{|z|}{{\mbox{\Small{\rm Re}}} z}%
\right)^{2\delta},
\end{equation*}
which is the claim.
}
\end{proof}

\begin{proof}[Proof of Theorem~$\protect\ref{gl}$]
\textrm{As before fix $x,y\in M$ and, for $0<2s< {d}(x,y)$, put $r={d}(x,y)-2s$%
. Next fix $f \in L^2(B(x,s), d\mu)\cap L^1(B(x,s), d\mu)$ and let $X\in TM$ be a smooth
vector field on $M$ supported in $B(y,s)$. This time  we set, for $z\in%
\mathbb{C}_+$,
\begin{equation*}  \label{ab1}
F(z)= \langle \nabla\exp(-z\Delta)\mathbf{m}_{V(\, \cdot \,,r)}f, X\rangle,
\end{equation*}
where ${V}(x,r)=\mu(B(x,r))$. 
Now
\begin{equation*}
F(z)= \langle \exp(-z\Delta) \mathbf{m}_{V(\, \cdot \,,r)}f, \nabla^* X\rangle
=\langle \exp(-z\Delta) \mathbf{m}_{V(\, \cdot \,,r)}f, \mbox{\rm div}\, X\rangle,
\end{equation*}
Therefore, since
$\Psi(z)=e^{-z\Delta}$ satisfies (\ref{DG1}) and
(\ref{DG2}),
$F$ satisfies (\ref{a11}) and (\ref{a22}) with 
\begin{equation*}
A=\|\mathbf{m}_{V(\,
\cdot \,,r)}f\|_{2} \|\mbox{\rm div}\, X\|_{2} \quad\mbox{and}\quad \gamma=r^2/4.
\end{equation*}
Note that $A$ is finite since $\mbox{\rm div}\, X $ is smooth and supported in
$B(y,s)$.
Now assume in addition that $\|f\|_1=\| |X| \|_1=1$ and let $z\in \mathcal{C%
}_{r^2/4}$. By Lemma~\ref%
{lemgln},
\begin{eqnarray*}
|F(z)| &\le&\|\nabla \exp(-z\Delta)\mathbf{m}_{V(\, \cdot \,,r)}\|_{1 \to \infty}
\\
&\le& C({\mbox{\Small{\rm Re}}} z)^{-a} \left(\frac{|z|}{{%
\mbox{\Small{\rm
Re}}} z}\right)^{2\delta}\left(\frac{r^2}{{\mbox{\Small{\rm Re}}} z}%
\right)^{\delta/2} =Cr^{-2a-4\delta}\left(\frac{r^2}{{%
\mbox{\Small{\rm Re}}} z}\right)^{(5\delta/2)+a}|z|^{2\delta}
\end{eqnarray*}
for all $z\in \mathcal{C}_{r^2/4}$. Thus $F$ satisfies the assumptions of
Proposition \ref{tw2} with 
\begin{equation*}
B=Cr^{-2a-4\delta}, \quad \nu=5\delta+2a, \quad \gamma=r^2/4,
  \quad \exp(g(z))=z^{-2\delta},\quad \mbox{and} \mbox{ any } \beta\in (0,1).
\end{equation*}
Therefore 
\begin{eqnarray*}
|F(z)| \le eCr^{-2a-4\delta} \left(\frac{r^2}{2|z|}%
\right)^{5\delta+2a }|z|^{2\delta} \exp\left(-{\mbox{\Small{\rm Re}}} \frac{r^2}{4z}\right)
\end{eqnarray*}
for all $z\in \mathcal{C}_{r^2/4}$. 
An obvious modification of Claim \ref{ker2} yields
\begin{eqnarray*}
|\nabla p_z(x,y)|V\left(y,r\right)& \le&\sup_{x^{\prime }\in B(x,s)}\sup_{y^{\prime }\in B(y,s)} |\nabla p_{z}(x^{\prime },y^{\prime })|V\left(y^{\prime },r\right)\\&=&\sup \left\{ |F(z)|\colon \,\Vert f\Vert _{L^{1}(B(x,s),d\mu)}=\Vert
X\Vert _{L^{1}(B(y,s),d\mu)}=1\right\} \\&\le& Cr^{-2a} \left(\frac{r^2}{4|z|}\right)^{3\delta+2a} \exp\left(- {\mbox{\Small{\rm Re}}}\frac{r^2}{4z}\right),
\end{eqnarray*}
and letting $s$ go to $%
0 $ we obtain 
\begin{equation}  \label{sil}
|\nabla p_z(x,y)|V\left(y,d(x,y)\right) \le C({d}^2(x,y))^{-a} \left(\frac{ {%
d}^2(x,y)}{4|z|}\right)^{3\delta+2a} \exp\left(- {\mbox{\Small{\rm Re}}}\frac{{%
d}^2(x,y)}{4z}\right)
\end{equation}
for all $z\in \mathcal{C}_{d^2(x,y)/4}$. If $0 \le t \le {d}^2(x,y)/4$ then $%
V(y,\sqrt t) \le V\left(y,\frac{d(x,y)}{2}\right)$, so by (\ref{sil}) 
\begin{equation*}
|\nabla p_t(x,y)|V(y,\sqrt t)\le C t^{-a} \left(\frac{ {d}^2(x,y)}{4t}%
\right)^{3\delta+2a} \exp\left(- \frac{{d}^2(x,y)}{4t}%
\right)
\end{equation*}
for all $0 \le t \le {d}^2(x,y)/4$. For $t \ge {d}^2(x,y)/4$, (\ref{rt1}) is immediate from
 \eqref{lemgl2n}.
}
\end{proof}

\textrm{\emph{Remark :}  Note  that (\ref{sil})  also yields
complex time estimates for the gradient of the heat kernel. }

\subsection{\textrm{Families of operators without semigroup property}}

\textrm{\label{nosemi} An important advantage of the technique which we
discuss here is that we do not have to assume that the family $%
\{\Psi(z)\colon \, z\in \mathbb{C}_+\}$ under consideration has the
semigroup property. Hence we are able to apply our results to families
operators which can be defined by:  $
\Psi(z)=g(z)\exp(-zL)$, where $g\colon \mathbb{C}_+ \to \mathbb{C}$ is an analytic function; 
$\Psi(z)=\exp(-zL_1) \exp(-zL_2)$;
$\Psi(z)=\exp(-zL_1) -\exp(-zL_2)$ or some more complex formulae. To be
more precise, let us come back to the general metric measure space setting,
and consider an analytic family of operators $\{\Psi(z)\colon \, z\in 
\mathbb{C}_+\}$ acting on $L^2(M,d\mu)$. Next assume that (see (\ref{ct})) 
\begin{equation}  \label{ctns}
\|\Psi(z)\|_{1 \to \infty} \le K({\mbox{\Small{\rm Re}}} z)^{-D/2}, \quad
\forall z\in \mathbb{C}_+.
\end{equation}
By \cite[Theorem~6, p.503]{DS} we can define the kernel $p_z^{\Psi}$ of the
operator $\Psi(z)$ in the same way as in (\ref{kernel1}) and again assume
that $p_z^{\Psi}$ is a continuous function on $M^2$. Now we can state the
following version of Theorem~\ref{truemain2}. }

\begin{theorem}
\textrm{\label{ns} Suppose that the family $\{\Psi(z)\colon \, z\in \mathbb{C%
}_+\}$ satisfies conditions \eqref{DG1}, \eqref{DG2} and \eqref{ctns}. Then
the kernel $p_z^{\Psi}$, if continuous, satisfies estimates \eqref{ge2}. }
\end{theorem}

\textrm{The proof of Theorem~\ref{ns} is the same as the proof of Theorem~%
\ref{truemain2}. }

\textrm{\medskip }

\textrm{For instance, let $L$ be a self-adjoint uniform elliptic second
order differential operator in divergence form with periodic coefficients
acting on $L^{2}(\mathbb{R}^{n}),$ and let $L^{o}$ be the corresponding
homogenized operator. Next let $p_{z}$ and $p_{z}^{o}$ be the corresponding
heat kernels. Then a straightforward modification of the argument from \cite%
{Zhi} shows the following so-called Berry-Esseen type estimate 
\begin{equation}
|p_{z}(x,y)-p_{z}^{o}(x,y)|\leq K\min \left\{1,\frac{|z|}{({\mbox{\Small{\rm
Re}}}z)^{3/2}}\right\} (\Small{\rm Re} z)^{-n/2}\label{perio}
\end{equation}%
for all $z\in \mathbb{C}_{+}$, $x,y\in \mathbb{R}^{n}$. The following
consequence of Theorem~\ref{ns} can be used to obtain Gaussian bounds for
the expression $|p_{z}(x,y)-p_{z}^{o}(x,y)|$. }

\begin{example}
\textrm{\label{pe} Suppose that $L$ and $L^0$ are two generators of analytic
semigroups on $L^2(\mathbb{R}^n)$. Next assume that $L$ and $L^0$ satisfy
conditions \eqref{DG1}, \eqref{DG2} with the distances $d$ and $d^0$.
Finally suppose that the corresponding heat kernels satisfy estimate \eqref%
{perio} and set $\tilde d(x,y)=\min\{d(x,y),d^0(x,y)\}$. Then 
\begin{equation*}
|p_z(x,y)-p_z^0(x,y)| \le \frac{eK\min \left\{1,\frac{|z|}{({\mbox{\Small{\rm Re}}}
z)^{3/2}}\right\}}{({\mbox{\Small{\rm Re}}} z)^{n/2}} \left(1+{\mbox{\Small{\rm Re}}} 
\frac{{\tilde{d}}^2(x,y)}{4z}\right)^{\frac{n+3}{2}} \exp\left(-{%
\mbox{\Small{\rm
Re}}} \frac{{\tilde{d}}^2(x,y)}{4z}\right)
\end{equation*}
for all $z\in \mathbb{C}_+$, $x,y\in\mathbb{R}^n$. }
\end{example}

\begin{proof}
\textrm{It easy to note that if for $z \in \mathbb{C}_+$ we put $%
\Psi(z)=\exp(-zL)-\exp(-zL^0)$ then $\Psi$ satisfies (\ref{DG1}) and (\ref%
{DG2}) with the distance $\tilde{d}$ (and with constant 2). Hence Example~%
\ref{pe} follows from Proposition~\ref{tw2}, with $\exp(g(z))=z$ and $\nu=n+3$. 
Note that with this choice of $g$,  \eqref{a10} is satisfied for any $\beta>0$, for some constants 
$C$ and $c$ depending on $\gamma$ and $\beta$.  }
\end{proof}

\textrm{\emph{Remark :} As we already said, the fact that the exponent in the
polynomial correction factor in front of the exponential in \eqref{ge2}
cannot be improved to $(D-1)/2$ as in \cite{Si1} is related to the fact that
the proof of Theorem~\ref{truemain2} does not use the fact that the family
of operators under consideration is a semigroup, in other words it is
related to the possibility of extending Theorem~\ref{truemain2} to Theorem~%
\ref{ns}. Indeed, let $L$ be the standard Laplace operator acting on $%
\mathbb{R}^n$, and apply Theorem~\ref{ns} to the family of operators $%
\Psi(z)=z^{n/2}\exp(-zL)$, $z\in \mathbb{C}_+$. In this case $D=0$ in %
\eqref{ctns}, and the conclusion cannot hold with $(D-1)/2=-1/2$ in the
polynomial correction factor in front of the exponential, since in that case 
$p_z(x,y)$ is exactly given by the Gauss function. By contrast, the argument
from \cite{Si1} cannot be applied to this choice of $\Psi(z)$, because it
only applies to semigroups. }

\textrm{A more elementary example is the following. Let $M=\{x,y\}$ with
counting measure and let ${d}(x,y)=1$. Consider the analytic family of
operators $\{ \Psi(z)\colon \, z \in \mathbb{C}_+\}$ acting on $L^2(M)$
given by the kernel $p_z(x,x)=p_z(y,y)=0$ and $p_z(x,y)=p_z(y,x)=\exp\left(-%
\frac{1}{4z}\right)$. It is easy to check that the family $\{\Psi(z)\colon
\, z \in \mathbb{C}_+\}$ satisfies conditions (\ref{DG1}),~(\ref{DG2})~and~(%
\ref{ctns}) with $D=0$. Again, this shows that the exponent $D/2$ cannot be
replaced by $(D-1)/2$ in the setting of Theorem~\ref{ns}. }

\subsection{\textrm{$L^p \to L^q$ Gaussian estimates}}

\textrm{Claim~\ref{ker1} reduces the proof of Gaussian bounds for the heat
kernel to obtaining a Gaussian type estimate for expressions of the form 
\begin{equation*}
\sup\left\{ |\langle \exp(-zL) f_1,f_2\rangle| \colon \, \|f_1\|_{L^1(U_1, d\mu)}
= \|f_2\|_{L^1(U_2, d\mu)}=1 \right\}.
\end{equation*}
In such expressions, one can replace the $L^1$ norms of functions $f_1$ and $%
f_2$ by the $L^p$ norm of $f_1$ and the $L^q$ norm of $f_2$, for $1\le p,q\le +\infty$. This leads to
natural generalizations of pointwise Gaussian bounds and provides some form
of Gaussian bounds for semigroups without heat kernels. There are many
interesting examples of operators which generate such semigroups. More
precisely, the corresponding semigroup $\exp(-tL)$ is not bounded from $L^1$
to $L^\infty$ even locally. The kernel of the operator $\exp(-tL)$ can always be
defined as a distribution or in some other sense, but in such cases it is not a
bounded function. Often such operators generate bounded semigroups on $L^p$
spaces only for $p$ ranging in some proper subinterval of $[1,\infty]$. We discuss a
semigroup of this type in Example~\ref{bez} below. }

\textrm{The so-called generalized Gaussian bounds that such semigroups may satisfy were  studied for instance  by Davies in
\cite[Lemmas~23~and 24]{Da1}, and they were extensively discussed  by Blunck and Kunstmann (see \cite{B1, B2, BK1, BK2, BK3}).
Estimates of a similar nature were also considered in \cite{LSV}, see
Propositions 2.6 and 2.8 in this paper.  For interesting considerations about $L^p-L^q$ Gaussian estimates, see also
 \cite{AM}.
}

\textrm{As in Section~\ref{nosemi}, we consider  analytic families of
operators $\{\Psi(z)\colon \, z\in \mathbb{C}_+\}$ on metric measure spaces
rather than semigroups generated by self-adjoint operators. Still the case $%
\Psi(z)=\exp(-zL)$ is the most natural example. In this section we are going
to consider families of operators satisfying the following condition 
\begin{equation}  \label{ondl}
\|\mathbf{m}_{{W_1}(\, \cdot \,, \sqrt{{\mbox{\Small{\rm Re}}} z})}\Psi(z) 
\mathbf{m}_{ {W_2}(\, \cdot \,, \sqrt{{\mbox{\Small{\rm Re}}} z})}\|_{p \to
q} \le 1, \qquad \forall z\in \mathbb{C}_+,
\end{equation}
where $1\le p<q\le +\infty$, and the functions $W_i$ satisfy condition~(\ref%
{d}) with exponents $\delta_i/2$ for $i=1,2$. We discuss the rationale for
condition~(\ref{ondl}) in remark (b) after Theorem~ \ref{th1}. Now let us
notice only that if $\Psi(z)=\exp(-zL)$, $p=1$, $q=\infty$ and $%
W_1=W_2=V^{1/2}$, then estimates~(\ref{ondl}) are equivalent to  estimates~(\ref{xy}), which follow as we have seen from condition~(%
\ref{on-dd}). Hence one can think of condition~(%
\ref{ondl}) as a generalization of the on-diagonal estimates (\ref{on-dd}). }

\begin{theorem}
\textrm{\label{th1} Let $(M,d,\mu)$ be a metric measure space. Suppose that
the functions $W_i\colon \mathbb{R}_+\times M \to \mathbb{R}_+$, $i=1,2$,
are continuous and satisfy condition \eqref{d} with constants $\sqrt{K'_i}\ge 1$ and exponents $%
\delta_i/2$. Next assume that the analytic family $\{\Psi(z) \colon \, z\in 
\mathbb{C}_+\}$ of operators on $L^2(M,d\mu)$ satisfies conditions \eqref{DG1}%
, \eqref{DG2} as well as \eqref{ondl} for some $p,q$, $1\le p<q\le +\infty$%
. Then 
\begin{equation}  \label{gau}
\|P_{1} \mathbf{m}_{{W}_1\left(\, \cdot \,,\frac{r}{2}\right)}\Psi(z)\mathbf{%
m}_{{W}_2\left(\, \cdot \,, \frac{r}{2}\right)} P_{2} \|_{p \to q}\le {e}
\sqrt{K'_1K'_2} \left(\frac{r^2}{4|z|}\right)^{\delta_1+\delta_2 }\exp\left(- {%
\mbox{\Small{\rm Re}}} \frac{r^2}{4z}\right)
\end{equation}
for all $U_1,U_2$ open subsets of $M$ and all $z\in \mathcal{C}_{r^2/4}$, with $r={d}%
(U_1,U_2)$. Here  $P_{i}$ denotes the operator of multiplication by the
characteristic function of the sets $U_i \subset M$, that is $P_i=
\mathbf{m}_{\chi_{U_i}}$ for $i=1,2$. 
 Moreover 
\begin{equation}  \label{gau1}
\|P_1 \mathbf{m}_{{W}_1(\, \cdot \,,\sqrt t)}\Psi(t)\mathbf{m}_{{W}_2(\,
\cdot \,, \sqrt t)} P_2 \|_{p \to q}\le {e} \sqrt{K'_1K'_2} \left(1+\frac{r^2}{4t}%
\right)^{\delta_1+\delta_2 } \exp{\left(-\frac{r^2}{4t}\right)}
\end{equation}
for all $t\in \mathbb{R}_+$. }
\end{theorem}

\begin{proof}
\textrm{Similarly as in (\ref{abnew}) we consider the function $F$ defined
by the formula 
\begin{equation*}
F(z)=\langle \Psi(z) \mathbf{m}_{W_2\left(\, \cdot \,,\frac{r}{2}%
\right)}f_2, \mathbf{m}_{W_1\left(\, \cdot \,,\frac{r}{2}\right)}f_1\rangle
\end{equation*}
but now assume that $f_i\in L^2(U_i,W_i^2(\, \cdot \,,r/2)\,d\mu)$ for $%
i=1,2 $ and $\|f_1\|_{L^{q^{\prime}}(U_1, d\mu)}= \|f_2\|_{L^{p}(U_2, d\mu)}=1$, where $1/q+1/q^{\prime}=1$%
. }

\textrm{In virtue of assumption~(\ref{DG1}) and Davies-Gaffney estimates (%
\ref{DG2}), $F$ satisfies (\ref{a11}) and (\ref{a22}) with $\gamma=r^2/4$
and 
\begin{equation*}
A= \| \mathbf{m}_{W_1\left(\, \cdot \,, \frac{r}{2}\right)}f_1\|_{2} \|%
\mathbf{m}_{W_2\left(\, \cdot \,, \frac{r}{2}\right)}f_2
\|_{2}=\|f_1\|_{L^2(U_1,W_1^2(\, \cdot \,,r/2)\,d\mu)} \|f_2
\|_{L^2(U_2,W_2^2(\, \cdot \,,r/2)\,d\mu)}.
\end{equation*}
Now if $z\in \mathcal{C}_{r^2/4}$, then $\sqrt{{%
\mbox{\Small{\rm Re}}} z}\le r/2$ by (\ref{real}). Using the assumptions on $f_1,f_2$, $%
W_1,W_2$ as well as (\ref{ondl}), we obtain 
\begin{eqnarray*}
|F(z)| &=& |\langle \mathbf{m}_{W_1\left(\, \cdot \,, \frac{r}{2}%
\right)}\Psi(z)\mathbf{m}_{W_2\left(\, \cdot \,, \frac{r}{2}\right)}f_2 ,
f_1 \rangle | \le \| \mathbf{m}_{W_1\left(\, \cdot \,, \frac{r}{2}\right)
}\Psi(z)\mathbf{m}_{W_2\left(\, \cdot \,, \frac{r}{2}\right)}\|_{p \to q} \\
&\le& \sup_{x,y\in M} \frac{W_1\left(x, \frac{r}{2}\right) W_2\left(y, \frac{%
r}{2}\right) }{{W}_1(x,\sqrt{{\mbox{\Small{\rm Re}}} z}) {W}_2(y,\sqrt{{%
\mbox{\Small{\rm Re}}} z}) } \| \mathbf{m}_{ {W}_1(\, \cdot \,,\sqrt{{%
\mbox{\Small{\rm Re}}} z})} \Psi(z)\mathbf{m}_{ {W}_2(\, \cdot \,,\sqrt{{%
\mbox{\Small{\rm Re}}} z}) }\|_{p \to q} \\
&\le& \sqrt{K'_1K'_2}\left(\frac{r^2}{{4{\mbox{\Small{\rm Re}}} z}}%
\right)^{(\delta_1+\delta_2)/{2}} =\sqrt{K'_1K'_2} 2^{-(\delta_1+\delta_2)} \left(%
\frac{{{\mbox{\Small{\rm Re}}} z}}{r^2}\right)^{-(\delta_1+\delta_2)/{2}}
\end{eqnarray*}
for all $z\in \mathcal{C}_{r^2/4}$. Thus $F$ satisfies (\ref{a33}) with $%
B=\sqrt{K'_1K'_2} 2^{-(\delta_1+\delta_2)}$ and $\nu=\delta_1+\delta_2$. By
Proposition~\ref{tw1}, 
\begin{equation*}
|F(z)| \le {e}\sqrt{K'_1K'_2} 2^{-(\delta_1+\delta_2)} \left(\frac{r^2}{2|z|}%
\right)^{\delta_1+\delta_2 }\exp\left(- {\mbox{\Small{\rm Re}}}\frac{r^2}{4z}%
\right).
\end{equation*}
The spaces $L^1(U_i, d\mu) \cap L^2(U_i,W^2_i(\, \cdot
\,,r/2)\,d\mu)$, $i=1,2$, being dense in $L^{q^{\prime}}(U_1, d\mu)$ (resp. $%
L^{p}(U_2, d\mu)$), the above inequality, for all functions $f_i\in
L^2(U_i,W^2_i(.,r/2)\,d\mu)$, $i=1,2$, such that $\|f_1\|_{q^{\prime}}=
\|f_2\|_{p}=1$, yields (\ref{gau}). }

\textrm{ To prove (\ref{gau1}) we notice that for $t\le r^2/4$, it is a
straightforward consequence of (\ref{gau}). For $t\ge r^2/4$, it follows
from (\ref{ondl}). }
\end{proof}

\textrm{\emph{Remarks :} (a) Note that if $\Psi(z)=\exp(-zL)$, $1/p+1/p'=1$,
and $W_1=W_2=W$, then it is enough to assume that (\ref{ondl}) holds for $%
z=t\in \mathbb{R}_+$. Indeed, by using, as in the proof of Theorem \ref{doudou}, the identity $\|T^*T\|_{p\to
p'}=\|T\|^2_{p\to 2}$ and the contractivity of $\exp(-isL)$, $s\in\mathbb{R}$%
, on $L^2(M,d\mu)$, one obtains 
\begin{equation*}
\|\mathbf{m}_{{W}(\, \cdot \,, \sqrt{{\mbox{\Small{\rm Re}}} z})}\exp(-zL) 
\mathbf{m}_{ {W}(\, \cdot \,, \sqrt{{\mbox{\Small{\rm Re}}} z})}\|_{p \to p'}
= \|\mathbf{m}_{{W}(\, \cdot \,, \sqrt{{\mbox{\Small{\rm Re}}} z})}\exp(- (%
{\mbox{\Small{\rm Re}}} z )L) \mathbf{m}_{ {W}(\, \cdot \,, \sqrt{{%
\mbox{\Small{\rm Re}}} z})}\|_{p \to p'}.
\end{equation*}
}

\textrm{(b) In \cite{B1, B2, BK1, BK2, BK3}, Blunck and Kunstmann develop
spectral multiplier theorems for operators which generate semigroups without
heat kernel acting on spaces satisfying the doubling condition\footnote{%
\textrm{Here we discuss only second-order operators. Blunck and Kunstmann
consider also the mth-order version of generalized Gaussian estimates.}}. As
their basic assumption they consider the following form of generalized
Gaussian estimates 
\begin{equation}  \label{bkgg}
\| \mathbf{m}_{\chi_{B(x,\sqrt{t})}}\exp(-tL) \mathbf{m}_{\chi_{B(y,\sqrt{t}%
)}} \|_{p \to q}\le C V(x,\sqrt t)^{\frac{1}{q}-\frac{1}{p}} \exp{\left(-c\frac{{d}^2(x,y)}{t%
}\right)},\ \forall\,t>0,\,x,y\in M,
\end{equation}
where $1\le p\le 2 \le q\le +\infty$ and $V(x,r)= \mu(B(x,r))$. The above
estimates imply that 
\begin{equation}  \label{bkond}
\| \exp(-tL) \mathbf{m}_{V^{(1/2)-(1/p)}(\, \cdot \,,
\sqrt{t})} \|_{p \to 2}\le
C \quad \mbox{and} \quad \| \mathbf{m}_{V^{(1/q)-(1/2)}(\, \cdot \,, \sqrt{t})}\exp(-tL) \|_{2 \to q}\le C
\end{equation}
(see \cite[Proposition~2.1, (ii)]{BK3}).
It follows from the above considerations that the opposite implication is
valid if $(M,d,\mu,L)$ satisfies conditions~(\ref{DG1}) and (\ref{DG2}).
Indeed, estimates (\ref{bkond}) imply (\ref{ondl}) with $W_1=V^{(1/2)-(1/q)}(\, \cdot \,, \sqrt{t})$ and $W_2=V^{(1/p)-(1/2)}(\,
\cdot \,, \sqrt{t})$. Now,
according to Theorem~\ref{th1}, (\ref{bkond}) implies the estimates (\ref%
{gau1}), which in turn imply (\ref{bkgg}) by choosing $U_1=B(x,\sqrt{t})$, $%
U_2=B(y,\sqrt{t})$ and using doubling. Thus Theorem~\ref{th1} can be used to
verify the main assumption of the results obtained in \cite{B1, B2, BK1, BK2, BK3}%
. For example the next statement follows from \cite[Theorem~1.1]{B2} and
Theorem~\ref{th1} (see also \cite[Theorem~4.3]{CCO}, for a more primitive version, with an additional $\varepsilon$ in the resulting exponent).  We give below  a proof that follows directly from  Theorem \ref{th1}.}  The conclusion of the corollary is instrumental in the theory of Riesz means (see \cite{CCO, B2}, and references therein);  see also condition $(HG_\alpha)$, p.339 in  \cite{GP} and its consequences.

\begin{coro}\label{ccob}
\textrm{Suppose that $(M,d,\mu,L)$ satisfies the doubling condition %
\eqref{dd}, as well as conditions \eqref{DG1}  and \eqref{DG2}. Let $p\in[1,2]$,  and assume that  there exists $C>0$ such that
\begin{equation}\label{jb}
\|\exp(-tL) \mathbf{m}_{V^{(1/p)-(1/2)}(\, \cdot \,, \sqrt{t})} \|_{p \to 2}\le
C, \quad \forall t\in \mathbb{R}_+,
\end{equation}
where $V(x,r)= \mu(B(x,r))$. Then there exists $C>0$ such that
\begin{equation*}
\|\exp(-zL) \|_{{\tilde p} \to {\tilde p}}\le C\left(\frac{|z|}{{{\mbox{\Small{\rm Re}}} z}}%
\right)^{\delta \left|\frac{1}{{\tilde p}}-\frac{1}{2}\right|}, \ \forall\,z \in \mathbb{C}_+, \, {\tilde p}\in
[p,p^{\prime}],
\end{equation*}
where $\delta$ is the exponent in condition \eqref{ddt}. }
\end{coro}

\begin{proof}
Note that,  for all $p\in[1,\infty]$, 
\begin{equation}  \label{ckl}
\sum_k \left(\sum_l |c_{lk} a_l|\right)^p \le \left(\max\left\{ \sup_l
\sum_k |c_{lk}|, \sup_k \sum_l |c_{lk}|\right\} \right)^p\sum_n |a_n|^p,
\end{equation}
with the obvious meaning for $p=\infty$, where $c_{lk}$, $a_l$ are sequences of real or complex numbers.
Indeed, for $p=1$ and $p=\infty$, (\ref{ckl}) is easy to obtain. Then we obtain (%
\ref{ckl}) for all $1\le p \le \infty$ by interpolation. Let $z\in \C_+$ and set $r^2=4({%
\mbox{\Small{\rm Re}}} z^{-1})^{-1}$, that is, $r=2\frac{|z|}{\sqrt{\mbox{\Small{\rm Re}}z}}$. Let $x_k$ be a maximal sequence in $M$ such
that all the balls $B(x_k,r/2)$ are disjoint. Note that the balls $B(x_k,r)$
are such that $\cup_k B(x_k,r)=M$, and that by doubling there exists $N\in\N^*$ such that any $x\in M$ is contained in at most $N$ such balls. Let $\chi_k$ be the characteristic
function of the set $B_k=B(x_k,r) \setminus \cup_{i=1}^{k-1} B(x_i,r) $, and set $r_{lk}=d(B_k,B_l)$. Using Jensen and doubling, we may write 
\begin{eqnarray*}
\|\exp(-zL)f\|_p^p&\le&\sum_k \left(\sum_l \|\chi_k \exp(-zL) (\chi_l f)\|_p\right)^p \\ &\le& C
\sum_k \left(\sum_l V(x_k,r)^{\frac{1}{p}-\frac{1}{2}} \|\mathbf{m}_{\chi_k} \exp(-zL) \mathbf{m}_{\chi_l}\|_{p \to 2}
\|\chi_l f\|_p\right)^p,
\end{eqnarray*}
hence by (\ref{ckl}) 
$$\|\exp(-zL)f\|_p^p\le C\max\{I,II\}^p \|f\|_p^p,$$
where 
\begin{eqnarray*}
I&=&\sup_l \sum_k V(x_k,r)^{\frac{1}{p}-\frac{1}{2}} \|\mathbf{m}_{\chi_k} \exp(-zL) \mathbf{m}_{\chi_l}\|_{p \to 2}, \\
II&=& \sup_k \sum_l V(x_k,r)^{\frac{1}{p}-\frac{1}{2}} \|\mathbf{m}_{\chi_k} \exp(-zL) \mathbf{m}_{\chi_l}\|_{p
\to 2}.
\end{eqnarray*}
Thus $\|\exp(-zL)f\|_{p\to p}\le C\max\{I,II\}$, and our goal is to estimate from above
$I$ and $II$ by $C\left(\frac{|z|}{{{\mbox{\Small{\rm Re}}} z}}%
\right)^{\delta \left|\frac{1}{p}-\frac{1}{2}\right|}$; the result for the other values of ${\tilde p}$ follows
 by duality and interpolation.
We shall explain how one deals with  $I$,  the treatment of $II$ being similar. 
Write
\begin{equation*}
\sup_l \sum_k V(x_k,r)^{\frac{1}{p}-\frac{1}{2}} \|\mathbf{m}_{\chi_k} \exp(-zL) \mathbf{m}_{\chi_l}\|_{p \to 2}=III+IV,
\end{equation*}
with
\begin{equation*}
III=\sup_l\sum_{\{k ; r_{kl}\le 2r\}} V(x_k,r)^{\frac{1}{p}-\frac{1}{2}} \|\mathbf{m}_{\chi_k} \exp(-zL) \mathbf{m}_{\chi_l}\|_{p
\to 2} 
\end{equation*}
and
\begin{equation*}
IV=\sup_l \sum_{i=1}^\infty\sum_{\{k ; 2ir < r_{kl} \le 2(i+1)r\}}V(x_k,r)^{\frac{1}{p}-\frac{1}{2}}
\|\mathbf{m}_{\chi_k} \exp(-zL)\mathbf{m}_{\chi_l}\|_{p \to 2}
\end{equation*}

Observe first that there exists $C$ only depending  on the doubling constant so that, for all $l$, $\#\{k; r_{kl}\le r\} \le C$.

Therefore
$$III\le 
C\sup_{\{k,l ;r_{kl} \le 2r\}}
 V(x_k,r)^{\frac{1}{p}-\frac{1}{2}} \|\mathbf{m}_{\chi_k}
\exp(-zL) \mathbf{m}_{\chi_l}\|_{p
\to 2}
$$ 

The conclusion of the corollary is trivial for $p=2$, therefore we can assume $1\le p<2$.  Using the contractivity of $\exp(-isL)$ on $L^2(M,d\mu)$,  one sees that \eqref{jb} implies \eqref{ondl},
 with $\Psi(z)=\exp(-zL)$, $q=2$, $W_1\equiv 1$ and $W_2=V^{(1/p)-(1/2)}$. 

Thus
\begin{eqnarray*}
 \|\mathbf{m}_{\chi_k}
\exp(-zL) \mathbf{m}_{\chi_l}\|_{p
\to 2}
&\le&\|\exp(-zL)\mathbf{m}_{\chi_l}\|_{p
\to 2}
\\
&\le&\|\exp(-zL)\mathbf{m}_{V^{(1/p)-(1/2)}(.,\sqrt{{\mbox{\Small{\rm Re}}} z})}\|_{p
\to 2}\|\mathbf{m}_{V^{(1/2)-(1/p)}(.,\sqrt{{\mbox{\Small{\rm Re}}} z})}\mathbf{m}_{\chi_l}\|_{p\to p}\
\\
&\le& \sup_{x\in B_l}
V^{\frac{1}{2}-\frac{1}{p}}(x,\sqrt{{\mbox{\Small{\rm Re}}}
z} ).
\end{eqnarray*}

Hence, if $r_{kl}\le 2r$,
\begin{eqnarray*}
III&\le& 
C\sup_{\{k,l ;r_{kl} \le 2r\}}
 V(x_k,r)^{\frac{1}{p}-\frac{1}{2}} \|\mathbf{m}_{\chi_k}
\exp(-zL) \mathbf{m}_{\chi_l}\|_{p
\to 2}
\\
&\le& C
\sup_{\{k,l ;r_{kl} \le 2r\}}\sup_{x\in B_l}
\left(\frac{V(x_k,r)}{V(x,\sqrt{{\mbox{\Small{\rm Re}}}
z} )}\right)^{\frac{1}{p}-\frac{1}{2}} \\
&\le& C
\sup_{\{k,l ;r_{kl} \le 2r\}}\sup_{x\in B_l}
\left(\frac{V(x,4r+r_{kl})}{V(x,\sqrt{{\mbox{\Small{\rm Re}}}
z} )}\right)^{\frac{1}{p}-\frac{1}{2}} 
\\ 
&\le& C
\left(\frac{V(x,6r)}{V(x,\sqrt{{\mbox{\Small{\rm Re}}}
z} )}\right)^{\frac{1}{p}-\frac{1}{2}} 
\\ &\le& C(r/\sqrt{{\mbox{\Small{\rm Re}}}
z})^{\delta\left(\frac{1}{p}-\frac{1}{2}\right)}\\
&\le& C(|z|/{\mbox{\Small{\rm Re}}}
z)^{\delta\left(\frac{1}{p}-\frac{1}{2}\right)}.
\end{eqnarray*}
In the last inequality, we could use doubling since $r=2\frac{|z|}{\sqrt{\mbox{\Small{\rm Re}}z}}\ge \sqrt{\mbox{\Small{\rm Re}}z} $. 

\medskip

Now for $IV$. Again, there exists $C$ only depending  on the doubling constant so that, for all $l$, $\#\{k; r_{kl}\le ir\} \le Ci^{\delta}$.

Therefore
$$IV\le C\sup_l\sum_{i=1}^\infty i^{\delta}\sup_{\{k ;2 ir< r_{kl} \le 2(i+1)r\}} V^{\frac{1}{p}-\frac{1}{2}} (x_k,r)\|\mathbf{m}_{\chi_k}
\exp(-zL)\mathbf{m}_{\chi_l}\|_{p \to 2} $$

Let us now estimate $\|\mathbf{m}_{\chi_k} \exp(-zL)\mathbf{m}_{\chi_l}\|_{p \to 2}$ for $r_{kl}>2r$ with the help of  Theorem \ref{th1}.
 Note that $\delta_1=0$ and $\delta_2=\delta\left(\frac{1}{p}-\frac{1}{2}\right)$, where $\delta$ is the exponent in \eqref{ddt}.

With our choice of $r$, $z\in \mathcal{C}_{r_{kl}^2/4}$ as soon as  $r_{kl}\ge r$.
Therefore  \eqref{gau} yields  \begin{eqnarray*}
&&\|\mathbf{m}_{\chi_k}
\exp(-zL)\mathbf{m}_{\chi_l}\|_{p \to 2} \\
&\le&\|\mathbf{m}_{\chi_k}
\exp(-zL)\mathbf{m}_{V^{(1/p)-(1/2)}(.,r_{kl}/2)}\mathbf{m}_{\chi_l}
\|_{p
\to 2}\|\mathbf{m}_{V^{(1/2)-(1/p)}(.,r_{kl}/2)}\mathbf{m}_{\chi_l}\|_{p\to p}
\\
&\le&C\left(\frac{r_{kl}^2}{4|z|}\right)^{\delta\left(\frac{1}{p}-\frac{1}{2}\right)}\exp\left(- {%
\mbox{\Small{\rm Re}}} \frac{r_{kl}^2}{4z}\right)\|\mathbf{m}_{V^{(1/2)-(1/p)}(.,r_{kl}/2)}\mathbf{m}_{\chi_l}\|_{p\to p}.
\end{eqnarray*}

Now
\begin{eqnarray*}
\|\mathbf{m}_{V^{(1/2)-(1/p)}(.,r_{kl}/2)}\mathbf{m}_{\chi_l}\|_{p\to p}&\le& \sup_{x\in B_l}V^{\frac{1}{2}-\frac{1}{p}}(x,r_{kl}/2)\\
&\le&( \inf_{x\in B_l}V(x,r_{kl}/2))^{\frac{1}{2}-\frac{1}{p}}\\
&\le& V^{\frac{1}{p}-\frac{1}{2}}(x_l,(r_{kl}/2)-r)\\
&\le& V^{\frac{1}{p}-\frac{1}{2}}(x_l,(i-1)r).
\end{eqnarray*}

Thus, if  $2 ir< r_{kl}\le 2(i+1)r$, 
\begin{eqnarray*}
\|\mathbf{m}_{\chi_k}
\exp(-zL)\mathbf{m}_{\chi_l}\|_{p \to 2} &\le&\left(\frac{r_{kl}^2}{4|z|}\right)^{\delta\left(\frac{1}{p}-\frac{1}{2}\right)}\exp\left(- {%
\mbox{\Small{\rm Re}}} \frac{r_{kl}^2}{4z}\right)V^{\frac{1}{2}-\frac{1}{p}}(x_l,(i-1)r)\\
&\le&V^{\frac{1}{2}-\frac{1}{p}}(x_l,(i-1)r)\left(\frac{(i+1)^2r^2}{|z|}\right)^{\delta\left(\frac{1}{p}-\frac{1}{2}\right)}\exp\left(- i^2r^2{%
\mbox{\Small{\rm Re}}} \frac{1}{z}\right)\\
&=&
 V^{\frac{1}{2}-\frac{1}{p}}(x_l,(i-1)r)\left(\frac{(i+1)^2r^2}{|z|}\right)^{\delta\left(\frac{1}{p}-\frac{1}{2}\right)}
e^{-4i^2} .
\end{eqnarray*}

Finally, 
\begin{equation*}
IV
\le C\sup_l\sum_{i=2}^\infty i^{\delta}\sup_{2 ir< r_{kl} \le 2(i+1)r}\left(\frac{V(x_k,r)}
{V(x_l,(i-1)r)}\right)^{\frac{1}{p}-\frac{1}{2}}\left(\frac{(i+1)^2r^2}{|z|}\right)^{\delta\left(\frac{1}{p}-\frac{1}{2}\right)}
e^{-4i^2}.
\end{equation*}

By doubling,
$$\left(\frac{V(x_k,r)}
{V(x_l,(i-1)r)}\right)^{\frac{1}{p}-\frac{1}{2}}\leq \left(\frac{V(x_l,3r+r_{kl})}
{V(x_l,(i-1)r)}\right)^{\frac{1}{p}-\frac{1}{2}}\leq \left(\frac{V(x_l,(2i+5)r)}
{V(x_l,(i-1)r)}\right)^{\frac{1}{p}-\frac{1}{2}}
$$
is uniformly bounded. Therefore
\begin{equation*}
 IV\le C \sum_{i=2}^\infty i^\delta \left(\frac{(i+1)^2r^2}{|z|}%
\right)^{\delta\left(\frac{1}{p}-\frac{1}{2}\right)} e^{-4i^2} \le C \left(\sum_{i=2}^\infty
e^{-4i^2/(1+\varepsilon)}\right) (|z|/{\mbox{\Small{\rm Re}}} z)^{\delta\left(\frac{1}{p}-\frac{1}{2}\right)},
\end{equation*}
which finishes   the proof.
\end{proof}

\medskip

We finish this section with the description of a simple and 
natural 
example of a family of operators which generate semigroups without heat 
kernel. We consider the following family of self-adjoint operators 
\begin{equation*}
L^{(c)}=\Delta-c|x|^{-2}
\end{equation*}
acting on $L^2(\mathbb{R}^n)$ for $n\ge 3$, where $\Delta=-\sum_{i=1}^n
\partial_{x_i}^2$ and $0,c \le (n-2)^2/4$.  Hardy's inequality shows that 
\begin{equation*} 
\Delta \ge (n-2)^2/4|x|^2 
\end{equation*} 
(see for example \cite[(2.1), p.107]{VZ}). Hence,  for all $c \in [0, (n-2)^2/4]$, $L^{(c)}$ is 
non-negative.
A detailed discussion of the definition of the operators  $L^{(c)}$
 can be found for example in \cite{AGG}.
  Such operators are
called Schr\"odinger operators with the inverse-square potential and they
are of substantial interest in analysis (see for example \cite{BPSTZ, VZ}
and references therein).
Note that $L^{(c)}$ is homogeneous of order $2$, meaning that if $U_t$ is the dilation
$(U_tf)(x)=f(tx)$, $t>0$, $x\in\R^n$, then  $U_{1/t}L^{(c)}U_t =t^2L^{(c)}$ for all $t>0$.
As a consequence,
$$
\|\exp(-tL^{(c)})\|_{q \to p}=t^{\frac{n}{2}\left(\frac{1}{q}-\frac{%
1}{p}\right)}\|\exp(-L^{(c)})\|_{q \to p}.
$$
Set 
$$
C_{c,p,q}=\|\exp(-L^{(c)})\|_{p \to q}.
$$
It was proved in \cite[Corollary~6.2]{VZ} that $C_{c,p,q}=\infty$ for all $p\le q$
and
$q>p^*_c =n/\sigma$, where $\sigma =(n-2)/2-\sqrt{(n-2)^2/4-c}$. It means in particular that for every $t>0$ the operator $%
\exp(-tL^{(c)})$ cannot be extended to a bounded operator on $L^p(%
\mathbb{R}^n)$ for $p>p*_c$, hence $L^{(c)}$ does not generate a 
semigroup on $L^p(\mathbb{R}^n)$ for such $p$. 
It also means  that for every $t>0$ the operator $\exp(-tL^{(c)})$ cannot
be extended to a bounded operator from $L^1(\mathbb{R}^n)$ to $L^\infty(%
\mathbb{R}^n)$: its kernel is not a bounded function on $%
\mathbb{R}^{2n}$ but merely a distribution. Therefore standard heat kernel
theory can not be applied to study the semigroup generated by $L^{(c)}$.
However, it was proved in  \cite{V}, see also  \cite{LSV},  that $\|\exp(-L^{(c)})\|_{p \to p}$ is finite for
all $p$ in the  interval $((p^*_ c)', p^*_c)$.

Our main interest here is some form of Gaussian type estimates which we
can
obtain for $L^{(c)}$ even though we know that that the pointwise Gaussian
estimates cannot hold. Our $L^p \to L^q$ Gaussian estimates are described
in the following example.
\begin{example} 
\textrm{\label{bez} Suppose that $c<(n-2)^2/4$ and that $%
p^{\prime}_* <p <2<q <p_*$, where $1/p_*+1/p^{\prime}_*=1$, 
$p_*=n/\sigma_c$ and $\sigma_c=(n-2)/2-\sqrt{(n-2)^2/4-c}$ . 
Then there exists a constant $C$ such that for any  two open subsets of 
$U_1,U_2 \subset \mathbb{R}^n$ the following estimates hold 
\begin{equation*} 
\| P_1 \exp(-zL^{(c)})P_2 \|_{p \to q}\le C({\mbox{\Small{\rm Re}}} 
z)^{-\nu} \left(1+{\mbox{\Small{\rm Re}}} \frac{r^2}{4z}\right)^{\nu} 
\exp\left(-{\mbox{\Small{\rm Re}}}\frac{r^2}{4z}\right),\ \forall\,z \in 
\mathbb{C}_+, 
\end{equation*} 
where $r={d}(U_1,U_2)$, ${d}$ is the Euclidean distance, 
$P_i=\mathbf{m}_{\chi_{U_i}}$ and 
$\nu=\frac{n}{2}\left(\frac{1}{p}-\frac{1}{q}\right)$. } 
\end{example} 

\begin{proof} 
\textrm{By Theorem~\ref{negpot}, $L^{(c)}$ 
satisfies Davies-Gaffney condition (\ref{DG2}), thus finite propagation 
speed for the corresponding wave equation with the standard Euclidean 
distance.  We need
\begin{equation} \label{dab}
\|\exp(-tL^{(c)})\|_{2 \to q}\le Ct^{\frac{n}{2}\left(\frac{1}{q}-\frac{%
1}{2}\right)}, \quad \forall t\in \mathbb{R}_+, 
\end{equation} 
for all $p$ such that $2\le p \le p_*$.  This is proved in 
\cite[Theorem~11~and~Lemma~13]{DaS}, but for the bounded potential  $1/(1+|x|^2)$ instead of
$1/|x|^2$. To circumvent this, we proceed as in  the proof of Theorem~\ref{negpot},  setting $L^{(c)}_a=-\Delta+%
\mathcal{V}_a$, where $\mathcal{V}_a=\max\{-c|x|^{-2},-a\}$, 
and we notice that, when $a$ goes to 
$%
\infty$, $L_a$ converges to $L+\mathcal{V}$ in the strong resolvent sense 
(see \cite[Theorem~VIII.3.3, p.454]{Ka} or \cite[Theorem S.16 p.373]{RS}). 
This yields \eqref{dab}. Moreover, 
\begin{equation*} 
\|\exp\left((-t+is)L^{(c)}\right)\|_{p \to q}\le 
\|\exp(-tL^{(c)}/2)\|_{2 \to q}\|\exp(-tL^{(c)}/2)\|_{p\to 2} 
\end{equation*} 
so 
\begin{equation} \label{skip}
\|\exp(-zL^{(c)})\|_{p \to q}\le C({\mbox{\Small{\rm Re}}} z)^{\frac{n}{%
2}\left(\frac{1}{q}-\frac{1}{p}\right)} 
\end{equation} 
for all $z\in \mathbb{C}_+$ and all $p,q$ such that $p^{\prime}_* <p \le q 
<p_*$. } 

\textrm{Hence for $p,q$ such that $p^{\prime}_* <q \le p <p_*$, $%
L^{(c)} $ satisfies all assumptions of Theorem~\ref{th1} with $%
W_1(x,t)=Ct^{\frac{n}{2}\left(\frac{1}{2}-\frac{1}{q}\right)}$ and $%
W_2(x,t)=Ct^{\frac{n}{2}\left(\frac{1}{p}-\frac{1}{2}\right)}$, 
$K'_1=K'_2=1$, 
$\delta_1=\frac{n}{2}\left(\frac{1}{2}-\frac{1}{q}\right)$ and $\delta_2=%
\frac{n}{2}\left(\frac{1}{p}-\frac{1}{2}\right)$. Thus Theorem~\ref{th1} 
yields 
\begin{eqnarray*} 
\| P_1 \exp(-zL^{(c)})P_2 \|_{p \to q}&\le& C r^{-(\delta_1+\delta_2)} %
\Big|\frac{r^2}{4z}\Big|^{\delta_1+\delta_2 }\exp\left(- {%
\mbox{\Small{\rm 
Re}}} \frac{r^2}{4z}\right) \\ 
&=& C({\mbox{\Small{\rm Re}}} z)^{-(\delta_1+\delta_2)} \Big({%
\mbox{\Small{\rm Re}}} \frac{r^2}{4z}\Big)^{(\delta_1+\delta_2)/2 
}\exp\left(- {\mbox{\Small{\rm Re}}} \frac{r^2}{4z}\right) 
\end{eqnarray*} 
for all $z\in \mathcal{C}_{r^2/4}$. The estimate for $z\not\in \mathcal{C}_{r^2/4}$ follows directly from \eqref{skip}. } 
\end{proof}

\subsection{\textrm{Other functional spaces}}

\textrm{\label{os} Our approach allows us to state and prove an analog of
off-diagonal bounds for other functional spaces. Suppose that $B$ is a
Banach space of functions on $M$ and that $L^2(M,d\mu)\cap B$ is dense in $B$ with respect to its
norm. Suppose next that $U$ is an open subset of $M$. We define the space $%
{B(}{U)}$ as the closure of $L^2(U, d\mu)\cap B$ in the space $B$.
Then we define ${B^*(U)}$ as the space $B^*$ with seminorm given by the
formula 
\begin{equation*}
\|f\|_{B^*(U)}=\sup_{\|g\|_{{B(}{U)}} \le 1}\langle f,g\rangle,
\end{equation*}
where $\langle f,g\rangle$  is the duality pairing.
}

\begin{theorem}
\textrm{\label{beso} Let $(M,d,\mu)$ be a metric measure space, $B_1$ and $%
B_2$ Banach spaces as above, and $U_1,U_2$ open subsets of $M$. Set $r={d}%
(U_1,U_2)$. Let $g\colon \mathbb{C}_+ \to \mathbb{C}$ be an analytic
function which satisfies condition \eqref{a10} for $\gamma=\frac{r^2}{4}$,
and let  $\{\Psi(z) \colon \, z\in \mathbb{C}_+\}$ be a family of bounded linear operators
 on $L^2(M,d\mu)$ satisfying
conditions \eqref{DG1}, \eqref{DG2}. Assume  that 
\begin{equation}  \label{bs}
\|\Psi(z)\|_{B_1 \to B^*_2} \le ({\mbox{\Small{\rm Re}}} z)^{-D/2}\exp%
\left(-{\mbox{\Small{\rm Re}}} g(z)\right), \ \forall \,z\in \mathbb{C}_+.
\end{equation}
Then, for all $z \in \mathbb{C}_+$, 
\begin{equation*}
\|\Psi(z)\|_{{B_1}{(U_1)} \to B^*_2(U_2)} \le ({%
\mbox{\Small{\rm
Re}}} z)^{-D/2} \left(1+{\mbox{\Small{\rm Re}}} \frac{r^2}{4z}%
\right)^{D/2}\exp\left(1-{\mbox{\Small{\rm Re}}} g(z)-{\mbox{\Small{\rm Re}}} \frac{r^2}{4z}\right).
\end{equation*}
}
\end{theorem}

About the growth condition \eqref {a10}, the same remark is in order than after Theorem \ref{grri}.

\begin{proof}
\textrm{For $f_i\in L^2(U_i, d\mu)\cap B_i$, $i=1,2$, we again consider the
function $F$ defined by the formula 
\begin{equation*}
F(z)=\langle \Psi(z) f_2, f_1\rangle.
\end{equation*}
By assumption 
\begin{equation*}
|F(z)|\le ({\mbox{\Small{\rm Re}}} z)^{-D/2}\exp\left({-{\mbox{\Small{\rm Re}}}
g(z)}\right) \|f_1\|_{{B_1}{(U_1)}} \|f_2\|_{{B_2}{(U_2)}}, \
z\in\mathbb{C}_+,
\end{equation*}
\begin{equation*}
|F(z)|\le \|f_1\|_{L^2(U_1, d\mu)}\|f_2\|_{L^2(U_2, d\mu)},\ z\in\mathbb{C}_+,
\end{equation*}
and 
\begin{equation*}
|F(t)|\le \exp\left(- \frac{r^2}{4t}\right)
\|f_1\|_{L^2(U_1, d\mu)}\|f_2\|_{L^2(U_2, d\mu)},\ t>0.
\end{equation*}
Therefore Proposition~\ref{tw2} yields 
\begin{equation*}
|F(z)|\le ({\mbox{\Small{\rm Re}}}
z)^{-D/2} \left(1+{\mbox{\Small{\rm Re}}} \frac{r^2}{4z}%
\right)^{D/2}\exp\left(-{\mbox{\Small{\rm Re}}} g(z)-{\mbox{\Small{\rm Re}}} \frac{r^2}{4z}\right)
\|f_1\|_{{B_1}{(U_1)}} \|f_2\|_{{B_2}{(U_2)}},
\end{equation*}
hence the claim. }
\end{proof}

\begin{example}
\textrm{\label{exa2} Let $\Delta$ be the Laplace-Beltrami operator acting on a
complete Riemannian manifold $M$. Assume that, for some $p\in(1,\infty)$ and
some $\alpha\ge 0$, 
\begin{equation}  \label{rel}
\|\exp(-z\Delta)\|_{p \to p} \le C\left(\frac{|z|}{{\mbox{\Small{\rm Re}}} z}%
\right)^{\alpha},\ \forall \,z\in \mathbb{C}_+.
\end{equation}
Then, for any pair $U_1,U_2$ of open subsets of $M$, one has 
\begin{equation*}
\left(\int_{U_2}|\exp(-z\Delta)f|^p\, d\mu\right)^{1/p} \le C \left(\frac{|z|}{{%
\mbox{\Small{\rm Re}}} z}\right)^{\alpha} \left(1+{\mbox{\Small{\rm Re}}} 
\frac{r^2}{4z}\right)^{\alpha} \exp\left(-{\mbox{\Small{\rm Re}}}\frac{r^2}{%
4z}\right) \|f\|^2_{L^p(U_1, d\mu)},
\end{equation*}
where $r={d}(U_1,U_2)$, for all $z\in \mathbb{C}_+$ and $f\in L^p(U_1, d\mu)$. }
\end{example}

\begin{proof}
\textrm{Note that in virtue of our assumptions the semigroup $\exp(-z\Delta)$, $%
z\in \mathbb{C}_+$ satisfies condition~(\ref{bs}) with $B_1=L^p(M)$, $B_2=L^{p'}(M) $,
$D=2\alpha$ and $|e^{g(z)}|=|z|^{\alpha}$. Hence Example~\ref{exa2}
follows from Theorem~\ref{beso}. }
\end{proof}

\textrm{For the relevance of assumption \eqref{rel}, see \cite{CCO}, \cite%
{B2}, and Corollary \ref{ccob} above. }

\bigskip

We can also set $g=\exp(\lambda z)$ in Theorem \ref{beso}. For $\lambda>0$, this allows  one to treat the case where one only has small time on-diagonal upper bounds,
for $\lambda>0$, this allows one, in case there is a spectral gap in the on-diagonal bounds, to keep track of it in the off-diagonal ones (see \cite{Ou} for more
in this direction). Such modification can be made in all our previous statements leading from on-diagonal bounds to off-diagonal bounds. We leave the details to the reader.

\subsection{\textrm{Possible further generalizations}}

\textrm{The technique presented above is very flexible. For instance, in the
statement and the proof of Theorem~\ref{th1} and many other results, we do
not have to assume that the operators $\Psi(z)$ are  linear. We think that it
should be possible to find    non-linear examples where the
Phragm\'en-Lindel\"of technique yields interesting results.  Also, other
pairs of dual norms than the $(L^2,L^2)$ norms could be considered in %
\eqref{DG1} and \eqref{DG2}. Of course, in such cases, one would loose the
connection with finite speed propagation. }

\subsection{\textrm{Open question}} Can one treat in a similar way the so-called sub-Gaussian estimates (see for instance \cite{BSF}, \cite{coupre}), which are typical of fractals, 
namely, can one imagine that, in the above notation,
$$
| \langle e^{-tL}f_1,f_2\rangle | \le \exp\left(-c\left(\frac{r^\beta}{t}\right)^{\frac{1}{\beta-1}}\right)
\|f_1\|_{2} \|f_2\|_{2},
$$
where $r$ is the distance between the supports of $f_1$ and $f_2$,
and 
$$
p_t(x,x)\le \frac{C}{V(x,t^{1\beta})},
$$
imply
$$
 p_t(x,y)\leq \frac{C'}{V(x,t^{1/\beta})}\exp\left(-c'\left(\frac{d^\beta(x,y)}{ t}\right)^{\frac{1}{\beta-1}}\right),
 $$
 for $\beta>2$?

\bigskip

{\bf Acknowledgements: } The first-named author would like to thank several institutions
where he could find the necessary peace of mind to work on this, in particular Macquarie University, Sydney, the University of Cyprus, the Isaac Newton Institute, Cambridge, and respectively  Xuan-Thinh Duong, Georgios Alexopoulos, and the organizers of the program on Spectral theory and partial differential equations, for giving him the opportunity to stay there.  The second-named author would like to thank Derek Robinson for organizing his 
extended visit to the Australian National University, where  part of this 
work was carried out. 
Both authors would like to thank Brian Davies, Alexander Grigor'yan, Alexander Teplyaev, whose questions helped to improve the manuscript.


\begin{thebibliography}{99}
\bibitem{A1} \textrm{G.~Alexopoulos. \newblock Spectral multipliers on {L}ie
groups of polynomial growth. \newblock {\em Proc. Amer. Math. Soc.},
120(3):973--979, 1994. }

\bibitem{A2} \textrm{G.~Alexopoulos. \newblock Oscillating multipliers on {L}%
ie groups and {R}iemannian manifolds. \newblock {\em Tohoku Math. J. (2)},
46(4):457--468, 1994. }

\bibitem{AGG}   \textrm{W.~Arendt, G.R.~Goldstein and J.A.~Goldstein. \newblock 
Outgrowths of Hardy's inequality, \newblock preprint, 2006.} 


\bibitem{ACDH} \textrm{P.~Auscher, T.~Coulhon, X.T.~Duong and S.~ Hofmann. %
\newblock Riesz transform on manifolds and heat kernel regularity. \newblock 
\emph{Ann. Sc. E. N. S.}, 37:911--957, 2004. }

\bibitem{AM} \textrm{P.~Auscher and J.-M.~Martell. \newblock Weighted norm inequalities, off-diagonal estimates and elliptic operators, Part II: Off-diagonal estimates on spaces of homogeneous type, \newblock preprint, 2005. }

\bibitem{BSF} \textrm{M.~Barlow. \newblock Diffusions on fractals, in {\em Lectures on probability theory and Statistics. Ecole
d'\'et\'e de probabilit\'es de St-Flour, XXV, 1995}, Springer Lecture Notes in Math., 1690, 1--121,
1998. }

\bibitem{BSC} \textrm{A. ~Bendikov and L.~Saloff-Coste.  \newblock On-
and off-diagonal heat kernel behaviors on certain infinite dimensional 
local Dirichlet spaces. \newblock {\em American J.  Math. }, 122,  1205--1263, 2000.}

\bibitem{B1} \textrm{S.~Blunck. \newblock A {H}\"ormander-type spectral
multiplier theorem for operators without heat kernel. 
\newblock {\em Ann.
Sc. Norm. Super. Pisa Cl. Sci. (5)}, 2(3):449--459, 2003. }

\bibitem{B2} \textrm{S.~Blunck. \newblock Generalized Gaussian estimates and
Riesz means of Schr\"odinger groups, \newblock unpublished manuscript, 2003. }

\bibitem{BK1} \textrm{S.~Blunck and P.C.~ Kunstmann. \newblock Calder\'on-{Z}%
ygmund theory for non-integral operators and the {$H\sp  \infty$} functional
calculus. \newblock {\em Rev. Mat. Iberoamericana}, 19(3):919--942, 2003. }

\bibitem{BK2} \textrm{S.~Blunck and P.C.~ Kunstmann. \newblock Weak type {$%
(p,p)$} estimates for {R}iesz transforms. \newblock {\em Math. Z.},
247(1):137--148, 2004. }

\bibitem{BK3} \textrm{S.~Blunck and P.C.~ Kunstmann. \newblock Generalized Gaussian estimates
and the Legendre transform. \newblock {\em J. Operator Th.},
53(2):351--365, 2005. }

\bibitem{BPSTZ} \textrm{N.~ Burq, F.~ Planchon, J.G.~ Stalker, and A.S.~
Tahvildar-Zadeh. \newblock Strichartz estimates for the wave and {S}%
chr\"odinger equations with the inverse-square potential. 
\newblock {\em J.
Funct. Anal.}, 203(2):519--549, 2003. }

\bibitem{CCO} \textrm{G.~Carron, T.~Coulhon, and E.-M.~Ouhabaz. \newblock %
Gaussian estimates and {$L\sp p$}-boundedness of {R}iesz means. \newblock    
\emph{J. Evol. Equ.}, 2(3):299--317, 2002. }

\bibitem{Ch} \textrm{M.~ Christ. \newblock {$L\sp p$} bounds for spectral
multipliers on nilpotent groups. \newblock {\em Trans. Amer. Math. Soc.},
328(1):73--81, 1991. }

\bibitem{Co1} \textrm{T. ~Coulhon. \newblock It\'eration de {M}oser et
estimation gaussienne du noyau de la chaleur. 
\newblock {\em J. Operator
Theory}, 29(1):157--165, 1993. }

\bibitem{coupre}
\textrm{T.~Coulhon.
\newblock Off-diagonal heat kernel lower bounds without Poincar\'e. 
\newblock{\em J. London Math. Soc.},  68(3):795--816, 2003.}


\bibitem{CD4} \textrm{T.~Coulhon and X.T.~ Duong. \newblock Riesz transforms
for {$1\leq p\leq 2$}. \newblock {\em Trans. Amer. Math. Soc.},
351(3):1151--1169, 1999. }

\bibitem{CD3} \textrm{T.~Coulhon and X. T.~Duong. \newblock Maximal
regularity and kernel bounds: observations on a theorem by {H}ieber and {P}%
r\"uss. \newblock {\em Adv. Diff. Eq.}, 5(1-3):343--368, 2000. }

\bibitem{CD2} \textrm{T.~Coulhon and X.T.~ Duong. \newblock Riesz transforms
for {$p>2$}. \newblock {\em C. R. Acad. Sci. Paris S\'er. I Math.},
332(11):975--980, 2001. }

\bibitem{CD1} \textrm{T.~Coulhon and X.T.~ Duong. \newblock Riesz transform
and related inequalities on noncompact {R}iemannian manifolds. \newblock    
\emph{Comm. Pure Appl. Math.}, 56(12):1728--1751, 2003. }

\bibitem{cgz} \textrm{T.~Coulhon, A.~Grigor'yan and F.~Zucca. The discrete
integral maximum principle and its applications. {\it Tohoku Math. J.}, 57(4):559--587, 2005. }

\bibitem{Da} \textrm{E. B.~ Davies. 
\newblock {\em Heat kernels and spectral
theory}. \newblock Cambridge University Press, Cambridge, 1989. }

\bibitem{Da2} \textrm{E. B.~Davies. \newblock Heat kernel bounds,
conservation of probability and the {F}eller property. 
\newblock {\em J.
Anal. Math.}, 58:99--119, 1992. \newblock Festschrift on the occasion of the
70th birthday of Shmuel Agmon. }

\bibitem{Da1} \textrm{E. B.~Davies. \newblock Uniformly elliptic operators
with measurable coefficients. \newblock {\em J. Funct. Anal.},
132(1):141--169, 1995. }

\bibitem{DP} \textrm{E. B.~Davies and M. M. H.~Pang. \newblock Sharp heat
kernel bounds for some {L}aplace operators. 
\newblock {\em Quart. J. Math.
Oxford Ser. (2)}, 40(159):281--290, 1989. }

\bibitem{DaS} \textrm{E. B.~Davies and B.~Simon. \newblock {$L\sp p$} norms
of noncritical {S}chr\"odinger semigroups. \newblock {\em J. Funct. Anal.},
102(1):95--115, 1991. }

\bibitem{DOS} \textrm{X. T.~Duong, E. M.~Ouhabaz, and A.~Sikora. \newblock %
Plancherel-type estimates and sharp spectral multipliers. 
\newblock {\em J.
Funct. Anal.}, 196(2):443--485, 2002. }


\bibitem{DS} \textrm{N.~ Dunford and J. T.~ Schwartz. 
\newblock {\em Linear
operators. {P}art {I}}. \newblock Wiley Classics Library. John Wiley \& Sons
Inc., New York, 1988. \newblock General theory, with the assistance of
William G. Bade and Robert G. Bartle, reprint of the 1958 original, a
Wiley-Interscience publication. }

\bibitem{Du} \textrm{N.~Dungey. \newblock Some remarks on gradient estimates for heat kernels. \newblock {\em Abstr. Appl. Anal.}, 73020,  2006. }


\bibitem{Ga} \textrm{M. P.~Gaffney. \newblock The conservation property of
the heat equation on {R}iemannian manifolds. 
\newblock {\em Comm. Pure Appl.
Math.}, 12:1--11, 1959. }

\bibitem{GP} \textrm{J.~Gal\'e and T.~Pytlik. \newblock Functional calculus for infinitesimal
generators of holomorphic semigroups. 
\newblock{\em  J. Funct. Anal.},  150(2): 307--355, 1997.}


\bibitem{Gr0} \textrm{A.~ Grigor{$^{\prime}$}yan. \newblock Heat kernel
upper bounds on a complete non-compact manifold. 
\newblock {\em Rev. Mat.
Iberoamericana}, 10(2):395--452, 1994. }

\bibitem{Gr1} \textrm{A.~Grigor{$^{\prime}$}yan. \newblock Gaussian upper
bounds for the heat kernel on arbitrary manifolds. 
\newblock {\em J. Diff.
Geom.}, 45(1):33--52, 1997. }

\bibitem{Gr} \textrm{A.~Grigor{$^{\prime}$}yan. \newblock Estimates of heat
kernels on {R}iemannian manifolds, \newblock in \emph{Spectral theory and
geometry (Edinburgh, 1998)},  London Math. Soc. Lecture
Note Ser.,  273, 140--225. Cambridge Univ. Press, Cambridge, 1999. }

\bibitem{Grf} \textrm{A.~Grigor{$^{\prime}$}yan. \newblock Heat kernel upper
bounds  on fractal spaces, 
\newblock preprint, 2004.  }


\bibitem{He} \textrm{W.~ Hebisch. \newblock Almost everywhere summability of
eigenfunction expansions associated to elliptic operators. 
\newblock {\em
Studia Math.}, 96(3):263--275, 1990. }

\bibitem{HR} \textrm{M.~Hino and J.~Ramirez. Small-time Gaussian behavior of
symmetric diffusion semigroups. \emph{Ann. Probab.}, 31(3):1254--1295,
2003. }

\bibitem{H} \textrm{L.~H\"ormander.   \newblock {\em
The analysis of linear partial differential operators. I. 
Distribution theory and Fourier analysis}. \newblock Second edition. Grundlehren der Mathematischen Wissenschaften 256. 
Springer-Verlag, Berlin, 1990. }


\bibitem{Ka} \textrm{T.~Kato. 
\newblock {\em Perturbation theory for linear
operators}. \newblock Classics in Mathematics. Springer-Verlag, Berlin,
1995. \newblock Reprint of the 1980 edition. }

\bibitem{KT} \textrm{M.~Keel and T.~Tao. \newblock Endpoint {S}trichartz
estimates. \newblock {\em Amer. J. Math.}, 120(5):955--980, 1998. }

\bibitem{LSV} \textrm{V.~Liskevich, Z.~Sobol, H.~Vogt. On $L_p$-theory of $%
C_0 $-semigroups associated with second order elliptic operators. \emph{J.
Funct. Anal.}, 193:55--76, 2002. }

\bibitem{MR} M.~Marias, E.~Russ. $H\sp 1$-boundedness of Riesz transforms and imaginary powers of the Laplacian on Riemannian manifolds. \emph{ Ark. Mat.}, 41(1):115--132, 2003. 

\bibitem{Ma} \textrm{A. I.~ Markushevich. 
\newblock {\em Theory of functions of a complex variable. {V}ol. {I}, {II},
  {III}}. \newblock Chelsea Publishing Co., New York, {E}nglish edition,
1977. \newblock Translated and edited by Richard A. Silverman. }

\bibitem{Mo} \textrm{S.~A. Molchanov. \newblock Diffusion processes, and {R}%
iemannian geometry. \newblock {\em Uspehi Mat. Nauk}, 30(1(181)):3--59,
1975. }

\bibitem{O} \textrm{E.~M. Ouhabaz.
\newblock {\em Analysis of heat equations on domains}. \newblock London Mathematical Society Monographs Series, 31. Princeton University Press, Princeton, NJ, 2005.
 }
 
 \bibitem{Ou} \textrm{E.~M. Ouhabaz.
\newblock Comportement des noyaux de la chaleur des op\'erateurs de Schr\"odinger et applications \`a certaines \'equations paraboliques semi-lin\'eaires, \newblock to appear in \emph{J.
Funct. Anal.} .
 }

   
\bibitem{RS} \textrm{M.~Reed and B.~ Simon. 
\newblock {\em Methods of modern
mathematical physics. {I}. Functional
analysis}. \newblock Academic Press Inc. [Harcourt Brace
Jovanovich Publishers], New York, second edition, 1980. }


\bibitem{Si1} \textrm{A.~Sikora. \newblock Sharp pointwise estimates on heat
kernels. \newblock {\em Quart. J. Math. Oxford Ser. (2)}, 47(187):371--382,
1996. }

\bibitem{Si3} \textrm{A.~Sikora. \newblock On-diagonal estimates on {S}%
chr\"odinger semigroup kernels and reduced heat kernels. 
\newblock {\em
Comm. Math. Phys.}, 188(1):233--249, 1997. }

\bibitem{Si2} \textrm{A.~Sikora. \newblock Riesz transform, Gaussian bounds
and the method of wave equation. \emph{Math. Z.}, 247(3):643-662, 2004. }

\bibitem{SW} \textrm{E.~M. Stein and G.~ Weiss. 
\newblock {\em Introduction
to {F}ourier analysis on {E}uclidean spaces}. \newblock Princeton University
Press, Princeton, N.J., 1971. \newblock Princeton Mathematical Series, No.
32. }

\bibitem{Stu1} \textrm{K.-Th.~Sturm. \newblock Analysis on local {D}irichlet
spaces. {II}. {U}pper {G}aussian estimates for the fundamental solutions of
parabolic equations. \newblock {\em Osaka J. Math.}, 32(2):275--312, 1995. }

\bibitem{Stu2} \textrm{K.-Th.~Sturm. \newblock The geometric aspect of {D}%
irichlet forms. \newblock In \emph{New directions in Dirichlet forms},
volume~8 of \emph{AMS/IP Stud. Adv. Math.}, pages 233--277. Amer. Math.
Soc., Providence, RI, 1998. }



\bibitem{VZ} \textrm{J. L.~Vazquez and E.~Zuazua. \newblock The {H}ardy
inequality and the asymptotic behaviour of the heat equation with an
inverse-square potential. \newblock {\em J. Funct. Anal.}, 173(1):103--153,
2000. }

\bibitem{V} \textrm{H.~Vogt. $L_p$-properties of second order elliptic 
differential equations, \newblock {\em PhD Thesis}, Dresden 2001.} 


\bibitem{Zhi} \textrm{V.V.~ Zhikov. Spectral approach to asymptotic problems
in diffusion. \newblock {\em Diff.\ Equations}, 25:33--39, 1989. }
\end{thebibliography}
\end{document}